\documentclass[11pt]{article}
\usepackage{fullpage}
\usepackage{palatino}
\usepackage{setspace,simplemargins}
\usepackage{latexsym,euscript,amsthm}
\usepackage{graphicx}
\usepackage{epsfig}
\usepackage{amsmath,amssymb,multicol,amscd}
\usepackage{mathrsfs}
\usepackage{framed}
\usepackage{subfig}
\newcommand{\bs}{\begin{eqnarray*}}
\newcommand{\es}{\end{eqnarray*}}

\newcommand{\bmat}{\begin{pmatrix} }
\newcommand{\emat}{\end{pmatrix}}

\newcommand{\beq}{\begin{equation}}
\newcommand{\eeq}{\end{equation}}
\newcommand{\beqa}{\begin{eqnarray}}

\newcommand{\eeqa}{\end{eqnarray}}

\newcommand{\integer}{\hbox{N \kern -.em I}\ }  
\newcommand{\real}{\hbox{R \kern -.2em I}\ }  
\newcommand{\numbersetmm}{\,\hbox{\scriptsize M  \kern -.2em I}\,\ }
\newcommand{\numbersetm}{\,\hbox{M \kern -.2em I}\,\,\ }

\newcommand{\bv}{{\mathbf v}}

\newcommand{\ep}{\epsilon}

\newcommand{\bx}{\mbox{\boldmath $x$}}
\newcommand{\field}[1]{\mathbb{#1}} 

\newcommand{\pr}{\partial}
\newcommand{\ba}{\mbox{\boldmath $a$}}

\newcommand{\bd}{\mbox{\boldmath $d$}}
\newcommand{\be}{\mbox{\boldmath $e$}}
\newcommand{\bU}{\mbox{\boldmath $U$}}
\newcommand{\bK}{\mbox{\boldmath $K$}}
\newcommand{\bT}{\mbox{\boldmath $T$}}

\newcommand{\bq}{\mbox{\boldmath $q$}}

\newcommand{\by}{\mbox{\boldmath $y$}}

\newcommand{\goto}{\rightarrow}
\newcommand{\B}{\mathbf}
\newtheorem{theorem}{Theorem}[section]

\newtheorem{definition}{Definition}[section]
\begin{document}

\begin{center}
{\fontsize{14}{20}\bf 
%
%
On reconstruction of dynamic permeability and tortuosity from data at  distinct frequencies%
}
\end{center}
\begin{center}
{\fontsize{12}{16} \bf 
%
Miao-Jung Yvonne Ou$^\dagger$, 
%
}
\end{center}
\begin{center}
{\fontsize{11}{12} \rm
%
$^\dagger${Department of Mathematical Sciences, University of Delaware, Newark, DE 19716, USA}, mou@math.udel.edu \\
}
\end{center}
 
{\fontsize{11}{12}\rm
\abstract
This article focuses on the mathematical problem of reconstructing the dynamic permeability $K(\omega)$ and dynamic tortuosity of poroelastic composites from permeability data at different frequencies, utilizing the analytic structure of the Stieltjes function representation of $K(\omega)$ derived by Avellaneda and Tortquato in \cite{avellaneda1991rigorous-link-b}, which is valid for all pore space geometry. The integral representation formula (IRF) for dynamic tortuosity is derived and its analytic structure exploited for reconstructing the function from a finite data set. All information of pore-space microstructure is contained in the measure of the IRF. The theory of multipoint Pad\'{e} approximates for Stieltjes functions guarantees the existence of relaxation kernels that can approximate the dynamic permeability function and the dynamic tortuosity function with high accuracy. In this paper, a numerical algorithm is proposed for computing the relaxation time and the corresponding strength for each element in the relaxation kernels. In the frequency domain, this approximation can be regarded as approximating the Stieltjes function by rational functions  with simple poles and positive residues. 
The main difference between this approach and the curve fitting approach is that the relaxation times and the strengths are computed from the partial fraction decomposition of the multipoint Pad\'{e} approximates, which is the main subject of the proposed approximation scheme.

With the idea from dehomogenization, we also established the exact relations between the moments of the positive measures in the IRFs of permeability and tortuosity  with  two important parameters in the theory of poroelasticity:  the infinite-frequency tortuosity $\alpha_\infty$ for the general case and the weighted volume-to-surface ratio $\Lambda$ for the JKD model, which is regarded as a special case of the general model. From these relations, we suggest a new way for evaluating these two microstructure-dependent parameters from a finite data set of permeability at different frequencies, without assuming any specific forms of the functions except the fact that they satisfies the IRFs. Numerical results for JKD permeability and tortuosity are presented. 
\vspace{0.2in}

\section{Notations}
The notations used in this paper are listed here.
\begin{itemize}
\item
 $A:=B$ means $A$ is defined by $B$
\item
$A=:B$ means $A$ defines $B$
\item
$A\approx B$ means $A$ is approximated by $B$
\item
Superscript $D$ is added to parameters to denote their values in the Johnson, Koplik and Dashen (JKD) model.
\item
$d\lambda$ is reserved for the positive measure in the IRF of permeability functions.
\item
$d\sigma$ is reserved for the positive measure in the IRF of tortuosity functions.
\end{itemize}
\section{\label{intro}Introduction}

Poroelasticity theory is a homogenized model for solid porous media
containing slightly compressible fluids that can flow through the pore structure.  This
field was pioneered by Maurice A. Biot, who developed his theory of
poroelasticity from the 1930s through the 1960s; a summary of much of
Biot's work can be found in his 1956 and 1962 papers ~\cite{biot1956theory-high,biot1956theory-of-propa,biot1962mechanics-of-de}.  Biot theory uses linear elasticity to describe
the solid portion of the medium (often termed the {\em skeleton} or
{\em matrix}), linearized compressible fluid dynamics to describe the
fluid portion, and Darcy's law to model the aggregate motion of the
fluid through the matrix.  While it was originally developed to model
fluid-saturated rock and soil, Biot theory has also been used in
underwater acoustics~\cite{bgwx:2004, gilbert-lin:1997,
  gilbert-ou:2003}, and to describe wave propagation in {\em in vivo}
bone~\cite{cowin1999bone-poroelasti, cowin-cardoso:fabric-2010,
  gilbert-guyenne-ou:bone-2012}. Biot's equations have been validated mathematically through applying homogenization techniques by various authors, e.g.  ~\cite{burridge1981poroelasticity-}~\cite{auriault1980dynamic-behavio,auriault1985dynamics-of-por} ~\cite{zhou1989first-principle}~\cite{gilbert2000homogenizing-th}~\cite{clopeau2001homogenizing-th}. Regularity of solutions to isotropic poroelastic equations has been studied in ~\cite{showalter2000diffusion}.

Biot theory predicts rich and complex wave phenomena within
poroelastic materials.  Three different types of waves appear: fast P
waves analogous to standard elastic P waves, in which the fluid and
matrix show little relative motion, and typically
compress or expand in phase with each other; shear waves analogous to
elastic S waves; and slow P waves, where the fluid expands while the
solid contracts, or vice versa.  The slow P waves exhibit substantial
relative motion between the solid and fluid compared to waves of the
other two types.  The viscosity of the fluid dissipates poroelastic waves as
they propagate through the medium, with the fast P and S waves being
lightly damped and the slow P wave strongly damped.  The viscous
dissipation also causes slight dispersion in the fast P and S waves,
and strong dispersion in the slow P wave.

A variety of different numerical approaches have been used to solve
poroelastic equations.  Carcione, Morency, and Santos provide a thorough
review of the previous
literature~\cite{carcione-morency-santos:comp-poro-rev}.  The earliest
numerical work in poroelasticity seems to be that of
Garg~\cite{garg-nayfeh-good:poro-1974}, using a finite difference
method in 1D.  Finite difference and pseudospectral methods have
continued to be popular since then, with further work by
Mikhailenko~\cite{mikhailenko:poro-1985},
Hassanzadeh~\cite{hassanzadeh:poro-1991}, Dai et
al.~\cite{dai-vafidis-kanasewich:poro-vel-stress}, and more recently
Chiavassa and Lombard~\cite{chiavassa-lombard:poro-fdm-2011}, among
others.  Finite element approaches began being used in the 1980s, with
Santos and Ore\~na's work~\cite{santos-orena:poro-fem} being one of
the first.  Boundary element methods have also been used, such as in
the work of Attenborough, Berry, and
Chen~\cite{attenborough-berry-chen:scattering}.  Spectral element
methods have also been used in both the frequency
domain~\cite{degrande-deroeck:poro-spectral-freq} and the time
domain~\cite{morency-tromp:poro-spectral-time}.  With the recent rise
of discontinuous Galerkin methods, DG has been applied to
poroelasticity in several works, such as that of de la Puente et
al.~\cite{delapuente-dumbser-kaser-igel:poro-dg}.  A Finite Volume Method solver for 2D and 3D Biot's equations can be found in \cite{lemoine2012high-resolution} and \cite{lemoine2013finite-volume-m}. There have also
been semi-analytical approaches to solving the poroelasticity
equations, such as that of Detournay and
Cheng~\cite{detournay-cheng:borehole}, who analytically obtain a
solution in the Laplace transform domain, but are forced to use an
approximate inversion procedure to return to the time domain.
Finally, there has been significant work on inverse problems in
poroelasticity, for which various forward solvers have been used ~\cite{sebaa2006ultrasonic-char}~\cite{buchanan2011recovery-of-the}; of
particular note is the paper of Buchanan, Gilbert, and
Khashanah~\cite{buchanan-gilbert-khashanah:bone-params-lofreq-2004},
who used the finite element method (specifically the FEMLAB software
package) to obtain time-harmonic solutions for cancellous bone as part
of an inversion scheme to estimate poroelastic material parameters,
and the later papers of Buchanan and
Gilbert~\cite{buchanan-gilbert:bone-params-hifreq-1-2007}, where the authors
instead used numerical contour integration of the Green's function. In \cite{sebaa2006ultrasonic-char}, numerical results from Biot's equations are compared with the experimental measurement of ultrasound propagation in cancellous bone. The physical parameters involved in the drag force are the tortuosity and the permeability. 

In the Biot equations for wave propagation in poroelastic materials \cite{biot1956theory-high,biot1956theory-of-propa}, a critical frequency $\omega_c$ is defined. For frequency below $\omega_c$, the pore fluid flow is laminar and  the friction term which takes into account the viscous interaction between the solid matrix and pore fluid is modeled by the product of friction constant and the difference between the fluid velocity and the solid velocity. We refer to this set of equations as low-frequency Biot equations. For frequency higher than $\omega_c$, the friction constant is multiplied by a frequency-dependent function to correct for the departure from laminar flow; this leads to a memory term in the time domain Biot equations. The exact form of the memory kernel is not known except for specific pore shapes such as parallel tubes \cite{biot1956theory-of-propa, auriault1985dynamics-of-por}. This set of equations are referred to as high-frequency Biot equations.

The need for quantifying the dissipation's dependance on frequency for more general pore space geometry prompted the work in the seminal paper \cite{johnson1987theory-of-dynam} by Johnson, Koplik and Dashen (JKD), in which the theory of dynamic fluid permeability $K^D(\omega)$ and  dynamic tortuosity $T^D(\omega)$ was developed for describing the {\bf inertial coupling} and {\bf viscous coupling} between matrix solid and pore fluid.  Using the causality argument, they derived the necessary symmetries and analytic properties of $K(\omega)$ and $T(\omega)$ when $\omega$ is extended to the complex values. Most importantly, they postulated the simplest forms of $K(\omega)$ and $T(\omega)$ which satisfy those properties. These two functions contain a tunable parameter $\Lambda$ take into account of the dependence on pore space geometries. However, it is very difficult to measure and is usually calculated through the empirical formula \cite{Pride1992Modeling-The-Dr}  $\Lambda\approx\sqrt{\frac{2\alpha_\infty K_0}{\phi/4}}$, where $\alpha_\infty$ is the infinite-frequency tortuosity, $K_0$ the static permeability and $\phi$ the porosity. The problem is that it is not clear how well this formula approximates $\Lambda$ and even if it does, the measurement of $\alpha_\infty$ is very difficult and is still an active research area \cite{allard1994evaluation-of-t}, \cite{dougherty2000a-quantitative-}, \cite{hughes2007investigation-o}, \cite{matyka2012how-to-calculat}, \cite{vallabh2010new-approach-fo}, \cite{dubeshter2014measurement-of-}.  
 
The Biot-JKD equations refer to the set of Biot equations modified by the JKD theory. In Biot-JKD equations no critical frequency is defined and the friction term (the drag force) is always a memory term.  

Due to the numerical complexity brought by the memory terms, most time-domain (vs frequency domain) solvers in the literature consider low-frequency Biot's equations even though it is well known in geological and biological applications that low-frequency Biot equations underestimate wave dissipation when compared with experiments. There have been a few papers which proposed different methods for handling the memory terms. Among them, the most popular ones are the fractional derivative approach for Biot-JKD model, which requires complicated quadrature rules \cite{lu2005wave-field-simu} and the phenomenological one which proposed to approximates the memory terms with sums of exponential decay kernels \cite{wilson1997simple-relaxati},\cite{carcione1996wave-propagatio}. The latter is more computationally efficient but it is not clear how the weights and decay rates of the exponential decay kernels can be found in a systematic way.

The aims of this paper are 
\begin{enumerate}
\item
To utilize the integral representation formula (IRF) of dynamic permeability, which is derived in \cite{avellaneda1991rigorous-link-b}, to develop a numerical scheme that can reconstruct the dynamic permeability function from any finite set of data measured at different frequencies. Unlike the JKD model, it does not impose any specific form on the permeability function. 
\item
To use the proposed numerical scheme, together with the relation between tortuosity and permeability, to reconstruct the tortuosity function from the finite data set. The weights and decay rates of the exponential decay kernel then come naturally along this process due to the mathematical structure of the tortuosity IRF derived in this paper.  
\item
To quantify how microstructure affects the tortuosity and other effective parameters relevant to drag force, which is known to be an important signaling mechanism for activating the cell process for bone remodeling,  \cite{robling2006biomechanical-a,johnson1996fluid-flow-stim, mcallister1999steady-and-tran,owan1997mechanotransduc,reich1991effect-of-flow-,sakai1998fluid-shear-str}. 
\end{enumerate}    

The paper is organized as follows. In Section \ref{background}, definition of permeability and its role in the poroelastic equations, the mathematical tools essential to the derivation of the tortuosity IRF and the inversion scheme are explained. In Section \ref{recon-perm}, the numerical scheme for reconstructing permeability functions from a finite data set is presented. Numerical results for the Biot-JKD model, which is regarded as a special case, are demonstrated. The proof that the JKD permeability function indeed can be represented as an IRF with a probability measure is also given there. In Section \ref{recon-tortu}, an IRF of the dynamic tortuosity function is derived. With this IRF, we prove that the time domain dynamic tortuosity function can be approximated by a combination of the Dirac function at $t=0$ and a sum of exponentially decay kernels whose rates and strengths can be computed from the proposed numerical scheme. Numerical results for the JKD tortuosity is demonstrated there. In Section \ref{moment-remark}, we present three exact (vs approximated) mathematical formulas which quantify how the geometry of pore space affects various effective poroelastic parameters through moments. Finally, in Section \ref{conclusion} we summarize the results and compare our exact formula for $\Lambda$ with an existing empirical formula. Also, future work is pointed out there. 
\section{\label{background}Mathematical Background}
\subsection{Permeability and Tortuosity}
For a {\bf rigid} porous medium filled with Newtonian pore fluid with density $\rho_f$ and {\bf dynamic} viscosity $\eta$, 
a key effective property is the fluid permeability tensor $\bK$, which is described by the so-called Darcy's law,
\cite{auriault1985dynamics-of-por,johnson1987theory-of-dynam,avellaneda1991rigorous-link-b,souzanchi2012tortuosity-and-} 
\beq
\bU = -\frac{\bK}{\eta} \nabla p
\eeq
where $\bU$ is the averaged fluid velocity over a representative volume element (RVE) of the porous medium and $\nabla p$ 
the applied pressure gradient; this is referred to as the {\bf static} permeability. If the applied pressure gradient is 
oscillatory with frequency $\omega$, then the induced averaged fluid velocity will also be oscillatory and proportional 
to $\nabla p(\omega)$ by
\beq
 \bU(\omega)= -\frac{\bK(\omega)}{\eta} \nabla p(\omega)
\eeq
where $\bK(\omega)$ is referred to as the {\bf dynamic} permeability, \cite{johnson1987theory-of-dynam,avellaneda1991rigorous-link-b}. 
$\bK(\omega)$ varies with $\omega$ because the viscous interaction between fluid and solid varies with frequency, 
as indicated by the frequency dependent viscous skin depth  $\sqrt{\frac{2 \eta}{\rho_f \omega}}$. 
For $\omega\ne 0$ the tortuosity tensor $\bT$ is related to $\bK(\omega)$ by 
\beq
\bT(\omega) =\frac{i \eta \phi}{\omega \rho_f} \bK^{-1}(\omega),\,\, i = \sqrt{-1} \mbox{  (note $\bT(\omega)$ has a pole at $\omega=0$).}
\label{tortuosity-permeability}
\eeq
As was mentioned in Section~\ref{intro},  there are two different forms of drag forces in Biot equations, depending on whether it is below 
or above the critical frequency $\omega_c:=\frac{\phi \eta}{\rho_f \alpha_\infty K_0}$, $K_0:=K(0)$, \cite{carcione2001wave-fields-in-}; 
for a low frequency, in the frequency domain, the drag force is $b\times$(fluid velocity relative to solid velocity) with $b=\eta \phi^2/K_0$, 
whereas for a high frequency, the constant $b$ is replaced by $b\cdot F(\omega)$ where $Im(F(\omega)) \goto 0$ and $Re(F(\omega)) \goto 1$ as $\omega\goto 0$. Biot derived the exact expression of $F(\omega)$ for thin circular tubes in terms of zero-order Kelvin functions of the first kind,
\cite{biot1956theory-high} and assumed the same functional form for all other pore geometry by a heuristically defined correction constant. 
In \cite{johnson1987theory-of-dynam}, the Biot-JKD equations are proposed by unifying the two types of friction terms in Biot's equations 
with a frequency-dependent function. For isotropic poroelastic materials, based on physics-based argument and exact calculation of parallel circular tubes, 
Johnson, Koplik and Dashen postulated the isotropic dynamic tortuosity to be of the form 
\beq
T^D(\omega)=\alpha_{\infty}\left(  1-\frac{\eta\phi}{i\omega\alpha_{\infty}\rho_{f} K_0}\sqrt{1-i\frac{4\alpha_{\infty}^{2}K_0^{2}\rho_{f}\omega}{\eta\Lambda^{2}\phi^{2}}}\right)=:{{\alpha(\omega)}}
\eeq
with the tunable geometry-dependent constant $\Lambda$, and (\ref{tortuosity-permeability}) implies  
\beq
K^D(\omega)= {K_0}/\left(\sqrt{1-\frac{4 i \alpha_\infty^2 K_0^2 \rho_f \omega}{\eta \Lambda^2 \phi^2}} -\frac{i \alpha_\infty K_0 \rho_f \omega}{\eta\phi}\right). 
\label{JKD-permeability}
\eeq
, where $\eta=\rho_f\,\nu$ is the dynamic viscosity of pore fluid.

The homogenization analysis in \cite{auriault1985dynamics-of-por} and physical arguments in \cite{johnson1987theory-of-dynam} shows that 
the permeability in Biot(-JKD) equations for poroelastic materials is identical to that for porous media with a rigid matrix. 
Furthermore, the permeability can be mathematically characterized as a functional of the solution to the unsteady Stokes equation,
 \cite{avellaneda1991rigorous-link-b}
\beq
\frac{\pr \bv}{\pr t} = -\nabla\left( \frac{p}{\rho_f}\right) + \nu \triangle \bv +v_0 \be \delta(t) \mbox{ in } \mathcal{V}_1,\,\
\nabla\cdot \bv = 0 \mbox{ in }\mathcal{V}_1\,\,
, \bv = {\bf 0} \mbox{ on }\pr \mathcal{V}
\label{LNSE-3}
\eeq 
where $\be$ is an arbitrary unit vector if $\bK$ is statistically isotropic, $v_0$ a constant, $\nu$ the {\bf kinetic} viscosity, $\delta(t)$ the Dirac delta function, 
$\mathcal{V}_1$ the region occupied by pore fluid and $\pr \mathcal{V}$ is the interface between fluid phase and solid phase in 
the RVE with periodic condition (or statistically homogeneous in the random media setting) on the outer boundary of RVE. As is indicated in \cite{avellaneda1991rigorous-link-b}, the results can be easily generalized to all statistically homogeneous anisotropic $\bK$. However, we assume $\bK$ is isotropic in this paper for simplicity.
It is shown in \cite{avellaneda1991rigorous-link-b} the solution $\bv(\bx,t)$ can be expressed as a sum of the normal modes $\Psi_n$
\beq
\bv(\bx,t) = v_0 \sum_{n=1}^\infty b_n e^{-t/\Theta_n} \Psi_n(\bx)
\label{stokes_expansion}
\eeq 
where $\Psi_n$ are the eigenfunctions of the Stokes system
\beq
\triangle \Psi_n +\nabla Q_n = -\ep_n \Psi_n \mbox{ and }\nabla\cdot \Psi_n =0 \mbox{ in }\mathcal{V}_1, \,\,\Psi_n =0 \mbox{ on } \pr \mathcal{V}_1, \Theta_n := (\nu \ep_n)^{-1} 
\eeq  
$0<\ep_1\le \ep_2\le \cdots$ and $\ep_n\goto \infty$ as $n\goto \infty$.  $\Theta_n$ are the viscous relaxation times and $\Theta_1$ referred to as the principal viscous relaxation time. The eigenfunctions are orthonormal in the sense
\beq
\frac{1}{|\mathcal{V}_1|} \int_{\mathcal{V}_1} \Psi_m(\bx) \cdot \Psi_n(\bx)\, d\bx = \delta_{mn} \mbox{ (Kronecker delta)}
\eeq
and the $b_n$ in (\ref{stokes_expansion}) 
\beq
b_n = \frac{1}{|\mathcal{V}_1|} \int_{\mathcal{V}_1} \be \cdot \Psi_n(\bx)\, d\bx 
\eeq
In \cite{avellaneda1991rigorous-link-b} it is shown through the classical Hodge decomposition argument that the infinite-frequency tortuosity 
$\alpha_\infty$ can be mathematically expressed as
\beq
{
\alpha_\infty = \phi F = \left(\sum_{n=1}^\infty b_n^2\right)^{-1}}
\eeq
This shows the microstructure information affects $\alpha_\infty$ through the projection of the applied flow direction on the normal modes 
of the Stokes equation in the pore space.
\subsection{Darcy's Laws in Poroelastic Equations and Permeability IRF}
The state variables for both 
Biot and Biot-JKD equations are $\bv$ (solid velocity), $\bq$ (fluid velocity relative to the solid) and $p$ (pore pressure). Note that $\bv$ and $p$ have different meanings from before and will stay unchanged hereafter. 
The stress-velocity formulation of Biot-JKD equations in a plane-strain case consists of 
\beqa
&& \pr_t \tau_{xx}=c_{11}^u \pr_x v_x+c_{13}^u\pr_z v_z+\alpha_1 M(\pr_x q_x+\pr_z q_z)+\pr_t s_1 \label{stress-strain-1} \\
&&\pr_t \tau_{zz}=c_{13}^u \pr_x v_x+c_{33}^u\pr_z v_z+\alpha_3 M(\pr_x q_x+\pr_z q_z)+\pr_t s_3\\
&&\pr_t \tau_{xz}=c_{55}^u (\pr_z v_x+\pr_x v_z)+\pr_t s_5\\
&&\pr_t p = -\alpha_1 M \pr_x v_x -\alpha_3 M \pr_z v_z -M(\pr_x q_x+\pr_z q_z)+\pr_t s_f,\label{stress-strain-4}\\
&&\rho \pr_t v_x + \rho_f \pr_t q_x = \pr_x \tau_{xx}+\pr_z\tau_{xz}\label{motion-x}\\
&&\rho \pr_t v_z + \rho_f \pr_t q_z = \pr_x \tau_{xz}+\pr_z\tau_{zz}\label{motion-z}
\eeqa
and the Darcy's laws which are the inverse Fourier transform of the following equation
\beq
-\nabla \tilde{p} = \eta \phi \bK^{-1} \tilde{\bq} +\rho_f (-i \omega) \tilde{\bv}_s = -i\omega\frac{ \rho_f}{\phi} \bT^D(\omega)\tilde{\bq}+\rho_f (-i \omega) \tilde{\bv}_s
\label{Darcy-freq}
\eeq 
where $\tilde{\bq}, \tilde{\bv}_s, \tilde{p}$ are Fourier transforms of $\bq$, $\bv$ and $p$, respectively. The Fourier transform we use here is
\[
\mathcal{F}[f](\omega)=\widetilde{f}(\omega):=\frac{1}{\sqrt{2\pi}}\int_0^\infty f(t) e^{i\omega t}\, dt
\]
The real part of $\bK^{-1}$ corresponds to the dissipation term and the imaginary part to the inertial term in the time domain where 
the drag force is expressed as a time-convolution term, \cite{auriault1985dynamics-of-por,lu2005wave-field-simu}. %
\begin{figure}[t]
\centering
\includegraphics[width=0.3\textwidth]{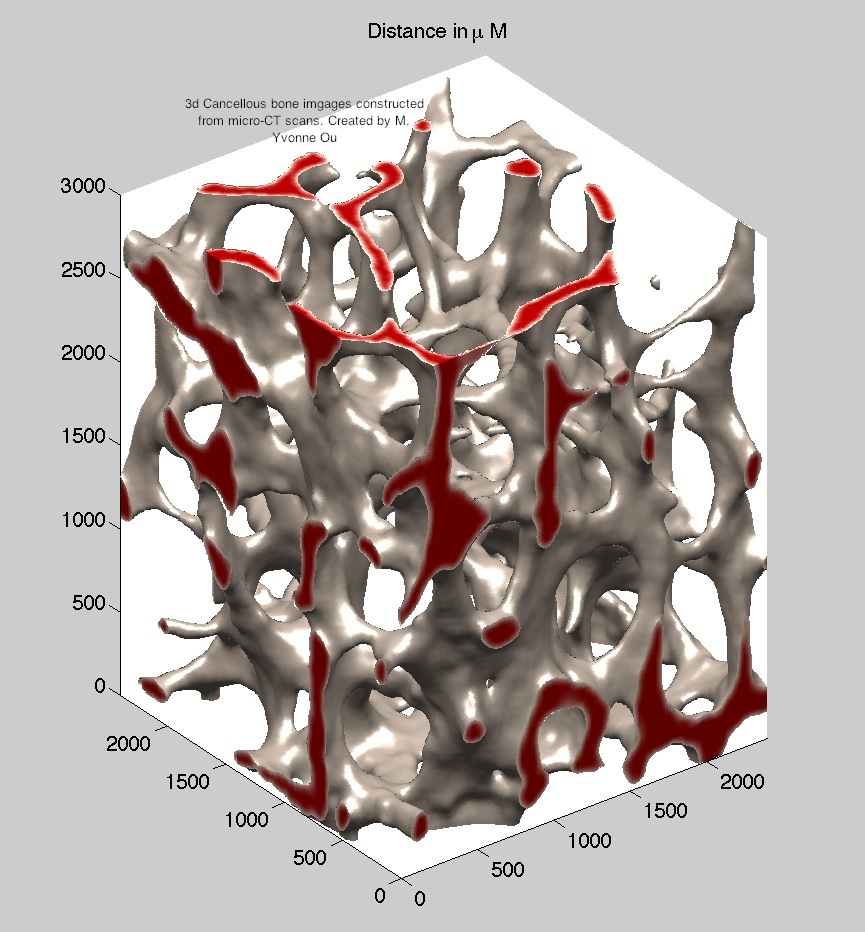}
\caption{\label{3Dbone}Cancellous bone, M.Y. Ou and L. Cardoso 2012}
\end{figure}
%
The complete form of any component of $\bK(\omega)$ for all $\omega$ is very difficult to compute for a given porous medium 
such as cancellous bone, whose pore geometry is complicated, see Fig.~\ref{3Dbone}. Hereafter, we consider isotropic $\bK$ and hence replace $\bK$ with $K$. It is pointed out in  \cite{johnson1987theory-of-dynam} that when $\alpha(\omega)$ and $K(\omega)$ are extended to the complex 
$\omega$-plane, they are analytic in the upper half plane because of causality and have the symmetry property 
$\alpha(-\overline{\omega})=\overline{\alpha(\omega)}$ and $K(-\overline{\omega})=\overline{K(\omega)}$, 
where the bar means complex conjugation. 
It is remarkable that in \cite{avellaneda1991rigorous-link-b} $K(\omega)$ is derived rigorously for {\bf all} pore space geometries 
as a Stieltjes integral with distribution $G(\Theta)$ that is nondecreasing, right-continuous, $G(\Theta)=0$ for $\Theta \le 0$ 
and $G(\Theta)=1$ for $\Theta \ge \Theta_1$ (i.e., a probability distribution) such that
\beq
 \left(\frac{F}{\nu} \right) K(\omega)=\int_0^{\infty} \frac{\Theta dG(\Theta)}{1-i\omega \Theta}, \mbox{ with } G(\Theta)=\frac{\sum_{\Theta_n \le \Theta} b_n^2}{\sum_{n=1}^\infty b_n^2}\,.
\label{K-omega-IRF}
\eeq
The derivation is based on identifying the exact functional form of $K(s)$ in terms of the solution to the Laplace transformed equations of  (\ref{LNSE-3}) 
with the parameter $s= - i\omega$.  
\subsection{Stieltjes functions and Multipoint Pad\'{e} approximation}
The reconstruction scheme proposed in this paper is based on the properties of $N$-point approximation for Stieltjes functions. Two different but related definitions of Stieltjes functions are widely used in the literature, we adopt the following definition in this paper.
\begin{definition}
A Stieltjes function $f(z)$ for $z$ on the extended complex plane has the following form
\beq
f(z)=\int_a^b \frac{d\mu(t)}{z-t}
\label{Stieltjes-function}
\eeq
where $a$, $b$ are extended real numbers and $\mu(t)$ is a bounded, non-decreasing real function.
\end{definition} 
A multipoint Pad\'{e} approximation of a function $f$ is a rational function interpolating $f$ at various points, not necessarily distinct. The following theorem (Theorem 1 on p. 26  of \cite{gelfgren1978rational}) is the foundation of the reconstruction algorithm proposed in this paper and hence we state it here. 
\begin{theorem}[\cite{gelfgren1978rational}]
Suppose $f$ is a Stieltjes function of the form in (\ref{Stieltjes-function}). Let $P_{n-1}(z)$ and $Q_n(z)$ be polynomials of degree at most $n-1$ and $n$, respectively, satisfying the relations ($k_1+k_2+k_3=2n$)
$$
\begin{cases}
f(z) Q_n(z)-P_{n-1}(z)= A(z) (z-x_1)\cdots(z-x_{k_1})(z-z_1)(z-\overline{z}_1)\cdots(z-z_{k_2})(z-\overline{z}_{k_2})\\
f(z) Q_n(z)-P_{n-1}(z)=B(z) z^{n-k_3-1}
\end{cases}
$$
where $A(z)$, $B(z)$ are analytic in $\field{C}\setminus[a,b]$, $B(z)$ bounded at $\infty$, $x_1,\cdots, x_{k_1}\in \field{R}\setminus [a,b]$, $z_1,\cdots, z_{k_2}\in \field{C}\setminus \field{R}$. Then $\frac{P_{n-1}(z)}{Q_n(z)}=\int_a^b \frac{d\beta(t)}{z-t}$ for some bounded, non-decreasing function $\beta(t)$. 
\label{N-pade-thm}
\end{theorem}
 Here $[n-1/n]_f(z):=\frac{P_{n-1}(z)}{Q_n(z)}$ is referred to as the $2n$-point Pad\'{e} approximant for $f$ and is unique~\cite{gelfgren1978rational}. It has nice convergent properties in $\field{C}\setminus[a,b]$ as long as the interpolating points are of the types specified in the theorem above. The significance of Theorem \ref{N-pade-thm} is that when the complex-valued interpolating  points appear in conjugate pairs, the approximant of a Stieltjes function can be expressed as
 \beq
 [n-1/n]_f(z)=\sum_{k=1}^{n} \frac{r_k}{z-p_k}
 \label{res-pole-R}
 \eeq
 with $r_k>0$ and $p_k\in(a,b)$ for $k=1,\cdots,n$. It will be made clear later in the paper that this property can be  used for generating efficient quadrature rules for dealing with the memory term in the dynamic Darcy's law (\ref{Darcy-freq}) in time domain.
\section{\label{recon-perm}Reconstruction of Dynamic Permeability}
Our experience with dehomogenization indicates that $K(\omega)$ can be reconstructed with a very good accuracy  from partial data by exploring its 
mathematical structure as a Stieltjes function, \cite{cherkaev2008dehomogenizatio,zhang2007inverse-electro,zhang2009reconstruction-}. 
Define a new variable $s$ and a new function $P(s)$ 
\beq
s:=-i\omega,\,\, P(s):= \left(\frac{F}{\nu} \right)K(is)=\int_0^{\Theta_1} \frac{\Theta dG(\Theta)}{1+s\Theta}=: \int_0^{\Theta_1} \frac{d\lambda(\Theta)}{1+s\Theta}
\label{P-function}
\eeq
\beq
\xi:=-1/s,\, R(\xi):=-s \left(\frac{F}{\nu} \right) K(is)= -sP(s)=\int_0^{\Theta_1} \frac{\Theta dG(\Theta)}{\xi-\Theta}\,.
\label{R-function}
\eeq
We summarize the definitions of auxiliary variables and auxiliary functions in Table \ref{transform}.
%
\begin{table}
\begin{center}
\begin{tabular}[t]{|c|c|c|}
\hline
$\omega$ & $s:=-i\omega$ & $\xi:=-\frac{1}{s}$ \\
\hline
$K(\omega)$ & $P(s):=(\frac{F}{\nu})\,K(is)$ & $R(\xi):=-sP(s)$\\
\hline
\end{tabular}
\end{center}
\caption{\label{transform} Auxiliary variables and functions for dynamic permeability $K(\omega)$}
\end{table}
%

It is clear that $R(\xi)$ in (\ref{R-function}) is a Stieltjes function with.
Due to the IRF of $P(s)$ in (\ref{P-function}), it is known that its Pad\'{e} approximants have accuracy-through-order property \cite{baker1996pade-approximan}, and all the poles of the Pad\'{e} approximants are simple with positive residue in $[-\infty,-\frac{1}{\Theta_1}]$ on the complex $s$-plane \cite{{ou2012on-nonstandard-}}. Most importantly, the Multipoint Pad\'{e} approximants (or rational interpolants) with interpolation knots $\{s_k\}_{k=1}^{N}$ with either $s_k >0$ or conjugate complex numbers appearing in pairs has interlacing simple zeros and simple poles locating in the regions where $P(s)$ is not analytic \cite{song2003convergence-of-}. 
Suppose we have values of $K(\omega_j)$ for different nonzero frequencies $\omega_1,\omega_2,\cdots, \omega_M \in \field{R}$. This means the values of $ R(\xi)$ for $\xi_k=\frac{-i}{\omega_k}\in \field{C}\setminus\field{R}$, $k=1,\cdots, M$ are known.  To generate the complex conjugated interpolating points at $\overline{\xi_k}$, we note that (\ref{R-function}) implies
\[
R(\overline{\xi})=\int_0^{\Theta_1} \frac{\theta dG(\theta)}{\overline{\xi}-\theta} = \overline{\int_0^{\Theta_1} \frac{\theta dG(\theta)}{\xi-\theta} } =\overline{R(\xi)}
\]
because $\theta,G(\theta)\in \field{R}$. Hence the data of $K$ at $M$ different frequencies indeed provide $2M$ data points for the reconstruction of $R$ through this symmetry. The $2M$-point Pad\'{e} approximants are formulated as follows.  
\beq
R(\xi_j) = \int_0^{\Theta_1} \frac{\Theta dG(\Theta)}{\xi_j-\Theta} \approx [M-1/M]_R(\xi):=\frac{a_0+a_1 \xi_j+\cdots+a_{M-1} \xi_j^{M-1}}{1+b_1 \xi_j+\cdots+b_{M} \xi_j^{M}},  j=1,\cdots, 2M
\label{LM-Pade}
\eeq
We know that the constant term in the denominator can be normalized to $1$ and the unknowns $a_0, \cdots, a_{M-1}, b_1,\cdots, b_M$ can be assumed real-valued because of (\ref{res-pole-R}).
Furthermore, the moments of $dG$ can be computed from partial fraction decomposition of it when lower frequency data points are used. Note that the first-moment of $dG$ is equal to  $\frac{F K_0}{\nu}$ and hence the formation factor $F:=\frac{\alpha_\infty}{\phi}$ can be recovered from the numerically estimated moments if $K_0$ and $\nu$ are known. That is, the tortuosity $\alpha_\infty$ can be recovered from data of $K(\omega)$ at low frequencies if the porosity $\phi$ is known.
In terms of the partial fraction decomposition of $[M-1/M]_R(\xi) $ 
\beq
[M-1/M]_R(\xi) = \sum_{j=1}^M \frac{r_j}{\xi-p_j},\, r_j>0, 0<p_j<\Theta_1
 \label{R-pol-res}
 \eeq 
the approximation of dynamic permeability $K$ can be expressed as
\beq
K(\omega) \approx \frac{\nu}{F} \sum_{j=1}^M \frac{r_j}{1-i\omega p_k}
\label{K-pol-res}
\eeq
Therefore, the permeability in time domain can be approximated as 
\beq
\mathcal{F}^{-1}[K](t)=\left( \frac{\nu}{F}\right)\sum_{j=1}^M \left(\frac{r_j}{p_j}\right) e^{-\frac{t}{p_j}}
\label{K-quad}
\eeq
\subsection{\label{R-reconstruct}Formulation and Algorithm}
For better conditioning of the inversion scheme, the reconstruction is based on (\ref{KD-IRF}), rather than $R(\xi)$.

Suppose we have values of $K(\omega_j)=K(i s_j)=P(s_j)$ for different nonzero  real-valued frequencies $\omega_1,\omega_2,\cdots, \omega_M$, then we can generate another $M$-interpolation points by using the symmetry of (\ref{K-omega-IRF}) for $\omega\in\field{R}$
\[
\left(\frac{F}{\nu} \right)  P(\overline{s_j})=K(-\omega_j)=K(-\overline{\omega_j})=\overline{K(\omega_j)}=\left(\frac{F}{\nu} \right)  \overline{P(s_j)}, j=1,\cdots, 2M.
\]
Because of (\ref{R-pol-res}) and (\ref{K-pol-res}), we can approximate $P(s)$ as 
\beq
P(s_j) = \int_0^{\Theta_1} \frac{\Theta dG(\Theta)}{1+s_j \Theta} \approx [M-1/M]_P(s):=\frac{a_0+a_1 s_j+\cdots+a_{M-1} s_j^{M-1}}{1+b_1 s_j+\cdots+b_{M} s_j^{M}}, \,\, j=1,2,\cdots, 2M
\label{PLM-Pade}
\eeq
and the moments of $dG$ can be computed from partial fraction decomposition of the approximant when lower frequency data are used.  We know that the constant term in the denominator can be normalized to $1$ because all the poles are simple and located in $(-\infty,-\frac{1}{\Theta_1})$. 

Given the data $(s_j,P(s_j))$, $s_j\ne 0$, $j=1,\cdots,M$, let $s_{j+M}=\overline{s_j}$, $P(s_{j+M})=\overline{P(s_j)}$ and $P_m:=P(s_m)$, $m=1,\cdots,2M$, then (\ref{PLM-Pade}) leads to the linear system of equations $A\bx=\bd$,
\beqa
A &=&
\begin{pmatrix}
1 & s_1 &s_1^2 &\cdots & s_1^{M-1} & -P_1 s_1 & - P_1 s_1^2 & -P_1 s_1^3 &\cdots & -P_1 s_1^{M}\\
1 & s_2 &s_2^2 &\cdots & s_2^{M-1} & -P_2 s_2  & -P_2 s_2^2 & -P_2 s_2^3 &\cdots & -P_2 s_2^{M}\\
\vdots &\vdots& \vdots&\vdots &\vdots& \vdots&\vdots &\vdots& \vdots&\vdots\\
1 & s_{2M} &s_{2M}^2 &\cdots & s_{2M}^{M-1} & -P_{2M} s_{2M} & -P_{2M}s_{2M}^2 & -P_{2M} s_{2M}^3 &\cdots & -P_{2M} s_{2M}^{M}
\end{pmatrix} \nonumber\\ &=& A_r+i\,A_i,\nonumber\\
\bx&=& (a_0, a_1,\cdots,a_{M-1},b_1, b_2, \cdots, b_M ), \nonumber\\
\bd&=&(P_1 \ P_2 \ \cdots\  P_{2M})^t
\label{linear-sys}
\eeqa
where $A_r$ and $A_i$ are the real part and imaginary part of $A$, respectively.
Since multi-point Pad\'{e} approximants of $P(s)$ constructed in this way have real-valued simple poles $p_j$ with positive residues $r_j$, $j=1,\cdots N$, we have 
\[
[M-1/M]_P(s):=\frac{a_0+a_1 s+\cdots+a_{M-1} s^{M-1}}{1+b_1 s+\cdots+b_{M} s^{M}}=\sum_{j=1}^M \frac{r_j}{s-p_j}
\]
and $\bx \in \field{R}^{2M}$.
 Noting that we have $2M$ real-valued unknowns and $2M$ complex-valued data for the linear system with complex coefficients, there are $4M$ equations with real-valued coefficients for the $2M$ real-valued unknowns. Rather than solving the formal equations of $A\bx=\bd$ as in \cite{zhang2007inverse-electro,zhang2009reconstruction-} in least-square sense, which requires forming $A_r^t A_r$ and $A_i^t A_i$, we solve this system of equations as an overdetermined least square problem by the following algorithm 
\begin{framed}
\begin{enumerate}
\item
Rescale each column $\ba_j$, $j=1,\cdots,2M$ of $A$ by $C:=diag\{ \|\ba_j\|_2^{-1} \}_{j=1}^{2M}$. Let 
\[A\bx=(AC)(C^{-1} \bx)=:B\by\]
\item
Let
\beq
\hat{B}:=\begin{pmatrix} B_r\\B_i \end{pmatrix}, \, \hat{\bd}:=\begin{pmatrix} \bd_r\\ \bd_i \end{pmatrix}.
\eeq

Solve the overdetermined system $\hat{B}\by=\hat{d}$ in the least square sense with Tikhonov regularization
\[
\min_{\by} \| \hat{B} \by - \hat{\bd} \|_2^2 +\gamma^2 \| \by \|_2^2 
\]
 The $L$-curve method \cite{hansen1993the-use-of-the-} is used for choosing the regularization parameter $\gamma$
\item
Rescale $\bx=C\by$.
\item
Apply partial fraction decomposition to the resulting $2M$-point Pad\'{e} approximant. Retain only the negative poles in $(-\infty, -\frac{1}{C_1}] \cup \{-\frac{1}{\xi_p} \}$ with positive weights and discard the rest.
\end{enumerate}
\end{framed}
Theoretically, all the poles in Step 4 should be in the range specified there. Numerically, the ill-posed nature of the inverse problem leads to poles outside the range. Suppose $M'$ poles are retained after discarding the spurious poles, $M' \le M $, we reindex them and the corresponding residues to $\{(p_j,r_j)\}|_{j=1}^{M'}$.  The function $P(s)$ is then approximated by
\beq
P(s)\approx P^M_{est}(s):=\sum_{j=1}^{M'} \frac{r_j}{s-p_j}=\sum_{j=1}^{M'} \frac{-r_j/p_j}{1+s(-1/p_j)}
\label{P-est}
\eeq
and the moments $\mu_k$, $k=0,1,2,\cdots$ by
\beq
\mu_k \approx  (-1)^{k+1} \sum_{j=1}^{M'}  \frac{r_j}{(p_j)^{k+1}}.
\label{momentE}
\eeq
In terms of the poles and residues in (\ref{P-est}), the time-domain permeability can be approximated by the relaxation kernel for $t \ge 0$
\[
\mathcal{F}^{-1}\{K\}(t)\approx \left(  \frac{\nu}{F}\right) \sum_{j=1}^{M'} r_j e^{p_j t}, \,\,r_j>0,\, p_j<0.
\]
Before testing the idea on the JKD permeability, we have to verify that it is consistent with the general theory presented in ~\cite{avellaneda1991rigorous-link-b}. 
\subsection{\label{IRF-JKD-Perm}IRF for JKD Permeability}
Since the JKD permeability in (\ref{JKD-permeability}) was derived by a completely different approach from that in ~\cite{avellaneda1991rigorous-link-b}, we need to show $K^D$ indeed assumes a representation of the form in (\ref{K-omega-IRF}). To see this, consider the auxiliary functions defined in Table \ref{transform} for $K^D$, 
\beq
R^D(\xi):=\frac{C_2}{\sqrt{\xi(\xi-C_1)}-C_2},\,\, C_2:=\frac{FK_0 }{\nu} >0,\,\,C_1:=\frac{4 C_2 F K_0}{\Lambda^2} >0,\,\, 
\eeq
and 
\[
P^D(s):= \left(\frac{F}{\nu} \right)K(is)=\left(\frac{F}{\nu} \right)\frac{K_0}{\sqrt{1+C_1 s}+C_2 s}
\]
%
For $R^D$ to assume the IRF, a specific branch of the square-root function has to be chosen so $R^D(\xi)$ has all the properties implied by the integral representation.  
The following branch for the square root function is chosen such that the branch cut of $R^D$ is contained in $[0,C_1]$, \cite{ablowitz2003complex-variabl}
\beq
R^D(\xi)=(r_1 r_2)^{1/2} e^{i(\theta_1+\theta_2)/2}
\eeq
where $(r_1,\theta_1)$ and $(r_2,\theta_2)$ are the local polar coordinates at the branch points $\xi=0$ and $\xi=C_1$ 
\beq
\xi=r_1 e^{i\theta_1},\,\, \xi-C_1=r_2 e^{i\theta_2},\,\, 0\le \theta_1, \theta_2 <2\pi
\eeq
With the chosen branch, the singular points of $R^D(\xi)$ consist of the branch cut $[0, C_1]$ and a simple pole at $\xi_p$, which is
\beq
\xi_p=\frac{C_1+\sqrt{C_1^2+4C_2^2}}{2}>C_1
\label{xip}
\eeq
\begin{figure}[t]
\centering
\includegraphics[scale=0.5]{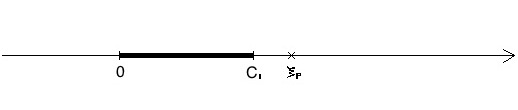}
\label{branch}
\caption{Branch cut and pole of $R^D$ on the complex $\xi$-plane}
\end{figure}
%
See Figure \ref{branch}. 

We would like to remark that 
\begin{enumerate}
\item
$-R^D(\xi)$ is analytic outside $[0,\xi_p]$
\item
$-R^D(\xi)$ maps the upper half plane to upper half plane with this choice of branch.
\item
There exists $\kappa>0$ such that $|y R^D(iy)|<\kappa$ for all $y>0$ because $|y R^D(iy)|\goto K_0$ as $\field{R}^+ \ni y\goto \infty$ and  $|y R^D(iy)|\goto 0$ as $\field{R}^+\ni y\goto 0$.  
\end{enumerate}
By a general representation theorem in function theory \cite{akhiezer1993theory-of-linea}, these three properties imply there exists a non-decreasing function $\lambda^D(u)$ of bounded variation on $[0,\xi_p]$ such that 
\beq
R^D(\xi)=\int_0^{\xi_p} \frac{d\lambda^D(u)}{\xi-u} \mbox{ for } \xi \in \field{C} \setminus [0,\xi_p]. 
\label{RD-IRF}
\eeq
Transforming (\ref{RD-IRF}) back to $K^D(i\omega)$, we see that $K^D$ indeed can be represented as an IRF with positive measure $d\lambda^D$
\beq
\frac{F}{\nu} K^D(\omega)=\int_0^{\xi_p} \frac{d\lambda^D(u)}{1-i\omega u}
\label{KD-IRF}
\eeq
Moreover, the IRF in (\ref{R-function}) gives additional information in the sense that it characterizes $d\lambda^D(u)$ as $u dG(u)$ with $dG$ being a probability measure. To verify this, we compute $d\lambda^D$ explicitly as follows. 

Since a function of bounded variation can only have jump discontinuities and is differentiable almost everywhere, $\lambda^D(u)$ must be continuous in the support of $\lambda^D(u)$ that corresponds to the branch cut of $R^D$. Using the Stieltjes inversion formula on page 224 of \cite{akhiezer1993theory-of-linea}, the density of $d\lambda^D$ in $[0,C_1]$ is
\beq
\psi(\xi):=\frac{d\lambda}{du} (\xi)=-\frac{1}{\pi} \lim_{y\goto 0^+} Im R(\xi+i y)=\frac{C_2 \sqrt{\xi (C_1 -\xi)} }{\pi\left[ C_2^2+\xi (C_1-\xi)\right]}
\label{psi}
\eeq
The pole of $R^D$ at $\xi=\xi_p$ corresponds to a Dirac measure of $d\lambda$ at $\xi_p$, whose strength $r$ can be computed from the following relation
\[
R^D(\xi)=\int_0^{C_1} \frac{\psi(u) du}{\xi-u} + \frac{r}{\xi-\xi_p}
\]
\beq
r=\lim_{\xi \goto \xi_p}\left[ R^D(\xi)  - \int_0^{C_1} \frac{\psi(u) du}{\xi-u} \right](\xi-\xi_p)=\lim_{\xi \goto \xi_p} R^D(\xi) (\xi-\xi_p) = \frac{2C_2 \sqrt{\xi_p(\xi_p-C_1)}}{2\xi_p-C_1}
\label{r}
\eeq
i.e. 
\[
d\lambda^D(u)=\chi_{\mathcal{I}}(u) \psi(u)  du+r\delta_{\xi_p},
\]
where $\chi_{\mathcal{I}}$ is the characteristic function of the interval $[0, C_1]$, $du$ the Lebesgue measure  and $\delta$ the Dirac measure. The following relation can be checked analytically
\beq
\int_0^{\xi_p} \frac{d\lambda^D(u)}{u} du = \int_0^{C_1} \frac{\psi(u)}{u}du+\frac{r}{\xi_p}=\frac{C_1}{\sqrt{C_1^2+4 C_2^2}}+\frac{4C_2^2}{\sqrt{C_1^2+4 C_2^2}(C_1+\sqrt{C_1^2+4 C_2^2})} =1
\eeq 
Our results of the JKD permeability IRF are summarized in the following theorem.
\begin{theorem}
\label{K-thm}
The JKD permeability in (\ref{JKD-permeability}) can be represented as
\[
K^D(\omega)=\frac{\nu}{F} \int_0^{\xi_p} \frac{u dG(u)}{1-i\omega u}
\]
where the probability measure $dG$ is
\[
dG(u)=\chi_{\mathcal{I}}(u) \left( \frac{\psi(u)}{u} \right)du+\left( \frac{r}{\xi_p} \right)\delta_{\xi_p},
\]
with $\xi_p$, $\psi$ and $r$ defined in (\ref{xip}), (\ref{psi}) and (\ref{r}), respectively.
\end{theorem}
\subsection{Numerical results for $K^D(\omega)$}
\label{NuRe-K}
We use the JKD permeability function given in (\ref{JKD-permeability}) to demonstrate the idea. 
\beq
P^D(s)=\frac{C_2}{C_2 s+\sqrt{1+C_1 s} } =\int_0^{\xi_p} \frac{d\lambda^D(u)}{1+su}
\label{PD-IRF}
\eeq
implies a specific branch of the square-root function has to be chosen so $P^D$ has all the properties implied by the integral representation such as
\begin{itemize}
\item
It maps $Im(s)>0$ to $Im(P(s))<0$.
\item
Its singularities are contained in $(-\infty,-\frac{1}{\Theta_1})$ for some $\Theta_1>0$.
\end{itemize}
We chose the branch with
\[
1+C_1 s=C_1(s-\left(-\frac{1}{C_1}\right)) = r e^{i\theta}, r=|1+C_1s| \mbox{ and } -\pi\le \theta <\pi
\]
and 
\[
\sqrt{1+C_1 s} =\sqrt{r} e^{i\theta/2}
\]
with branch cut at $(-\infty,-\frac{1}{C_1})$. 
With this choice of branch, it can be verified that
\begin{enumerate}
\item
$P^D(s)$ is analytic outside $(-\infty, -\frac{1}{C_1}] \cup \{\frac{C_1-\sqrt{C_1^2+4C_2^2}}{2C_2^2}=\frac{-1}{\xi_p}\}$.
\item
$P^D(s)$ maps the upper half plane to lower half plane.
 \end{enumerate}
\begin{table}[b]
\caption{\label{table-para} Parameters used in numerical simulations.}
\begin{center}
\begin{tabular}{|c|c|c|c|c|c|c|c|}
\hline
$\phi$ &  $\alpha_\infty$ & $K_0$$(m^2)$&$\nu=\eta/\rho_f$ & $\Lambda$ $(m)$ &$F$&$C_2$ &$C_1$\\
\hline
0.67 & 1.08 & $7\times 10^{-9}$ & $\frac{30\times10^{-3} Kg/(m\cdot s)}{1060 Kg/m^3}$ &  $10^{-5}$ & $\alpha_\infty/\phi$& $\frac{FK_0}{\nu}$ & $\frac{4C_2 F K_0}{\Lambda^2}$\\
\hline
\end{tabular}
\end{center}
\end{table}
The exact values of the moments of $d \lambda^D$ can be computed by differentiating (\ref{PD-IRF}) near $s=0$ and equating the coefficients on both sides.
\beq
\mu_k:=\int_0^{\xi_p} u^k d\lambda^D(u)=(-1)^k \frac{P^{(k)}(0)}{k!}
\label{momentNP}
\eeq

The JKD-Biot parameters of cancellous bone taken from literature \cite{cowin1999bone-poroelasti, hosokawa1997ultrasonic-wave, fellah2004ultrasonic-wave}: $\phi=0.67,\, \alpha_\infty=1.08,\,K_0=7\times 10^{-9}$ $m^2$, $\nu=\eta/\rho_f= \frac{30\times10^{-3} Kg/(m\cdot s)}{1060 Kg/m^3}$, $\Lambda=10^{-5}$ $m$, $F=\alpha_\infty/\phi$. This corresponds to $C_1=0.17994, C_2=3.98687\times 10^{-4}$, $-\frac{1}{C_1}=-5.55726$ and $-\frac{1}{\xi_p}=-5.55724$.

We first demonstrate the results in lower frequency range, when the reconstruction of moments is of interest. For $\omega$ between 1 Hz and 51 Hz, let $M=4,6,8,10$. For each fixed $M$, the frequency range $[1,51]$ are equally divided into $M-1$ intervals with sample  frequency taken at $\omega_0=1$, $\omega_2=1+\delta \omega, \cdots, \omega_M=51$ with $\delta\omega=50/(M-1)$. The corresponding moments computed from (\ref{momentE}), together with the exact moments are listed in Table \ref{table-moments}. The plot of these moments estimated with various values of $M$ is in Figure \ref{fig-table3}.

Figure {\ref{rel-err-1-51}} shows the max. relative error of estimating $P(s)$ with $P_M^{est}$ in (\ref{P-est}), where the maximum is taken among $1000$ equally spaced sample points. The maximum relative error $E_\infty$ is defined as follows.
\[
{\mbox{$E_\infty$($P^M_{est}$)}} := \max_{s\in [-i\omega_{min}, -i\omega_{max}]}  \frac{| P(s)-P^M_{est}(s)|}{|P(s)|}
\]
\begin{table}
\caption{\label{table-moments}Moments of $d\lambda^D$ from data of $\omega\in [1, 51]$ Hz.}
{\small
\begin{center}
\begin{tabular}{|c|c|c|c|c|c|}
\hline
& Exact & $M=10$ &  $M=8$ &  $M=6$ &   $M=4$ \\
\hline
$M'$ & & 9 & 7 &5 & 3\\
\hline
$\mu_0$ & 0.3986866 e-3& 0.3986865 e-3 &  0.3986831 e-3 & 0.3985660 e-3& 0.3951767 e-3\\
$\mu_1 $ & 0.3602968 e-4 & 0.3602966 e-4  &  0.3602864 e-4  & 0.3598232 e-4  & 0.3469776 e-4 \\
$\mu_2 $ & 0.4869721 e-5& 0.4869696 e-5&  0.4868684 e-5& 0.4837053 e-5 & 0.4286249 e-5\\
$\mu_3 $ & 0.7310984 e-6 & 0.7310788 e-6 &  0.7304975 e-6& 0.7173481 e-6& 0.5624193 e-6\\
$\mu_4 $ & 0.1152294. e-6 & 0.1152184. e-6&  0.1149672 e-6& 0.1107076 e-6& 0.7486713 e-7\\
$\mu_5 $ & 0.1867826 e-7 & 0.1867329 e-7& 0.1858360 e-7& 0.1740639 e-7 & 0.1000334 e-7\\
$\mu_6 $ & 0.3083493 e-8  & 0.3081572 e-8& 0.3053660 e-8& 0.2761889 e-8& 0.1337903 e-8\\
$\mu_7 $ & 0.5156154  e-9 & 0.5149587 e-9 & 0.5071242 e-9& 0.4402349 e-9 & 0.1789848 e-9 \\
$\mu_8 $ & 0.8704482 e-10 & 0.8684085 e-10 & 0.8481019 e-10 & 0.7033276 e-10 & 0.2394621 e-10 \\
$\mu_9 $ & 0.1480291 e-10 & 0.1474427 e-10 & 0.1424991 e-10&  0.1124946 e-10 &0.3203798 e-11\\
\hline
\end{tabular}
\end{center}
}
\end{table}  
\begin{figure}[h]
\centering
\includegraphics[width=0.5\textwidth]{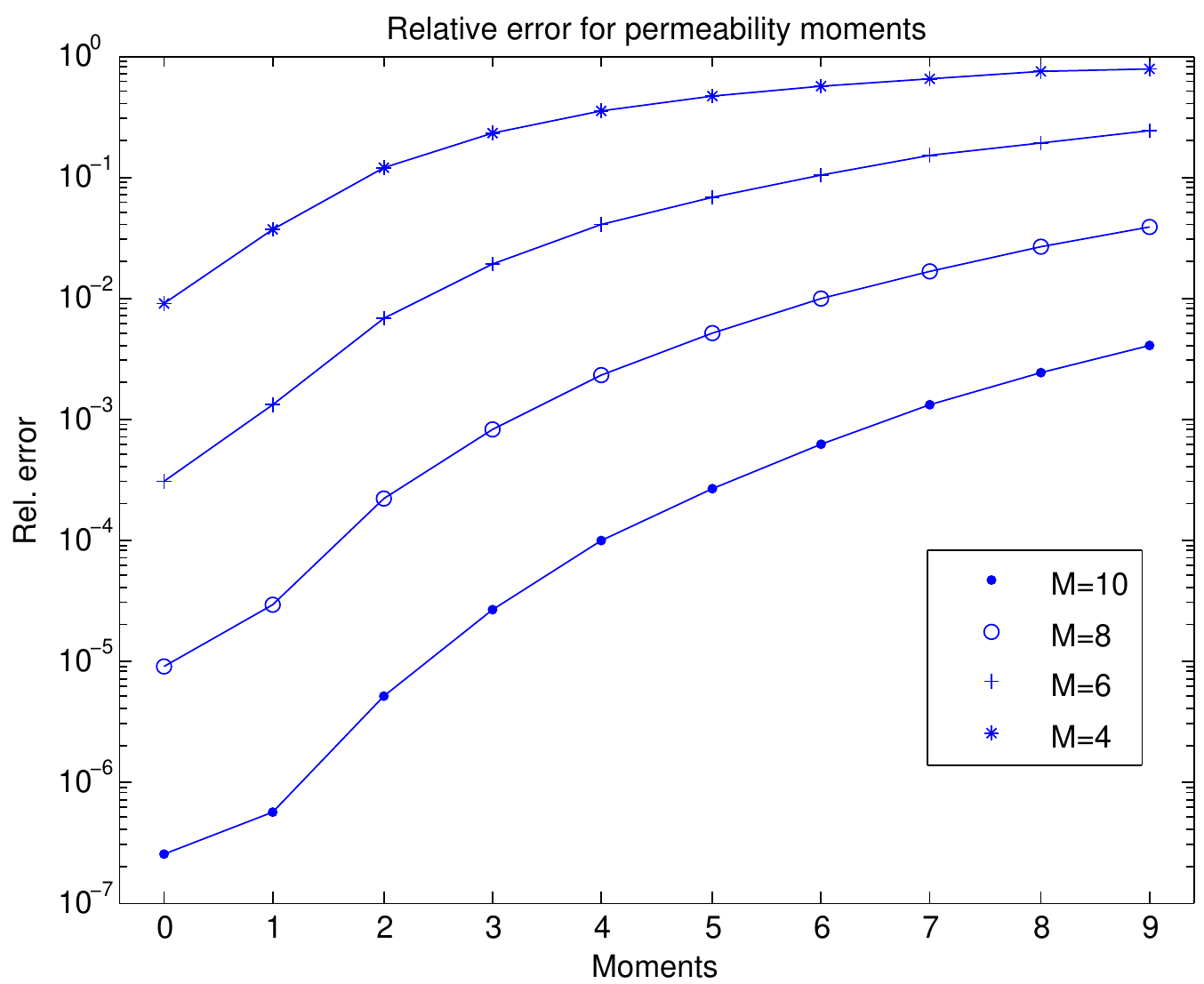}
\caption{\label{fig-table3} Relative error of $\mu_1(d\lambda^D), \cdots, \mu_9(d\lambda^D)$ approximated from $M$ data points}
\end{figure}
\begin{figure}[h]
\centering
\includegraphics[width=0.5\textwidth]{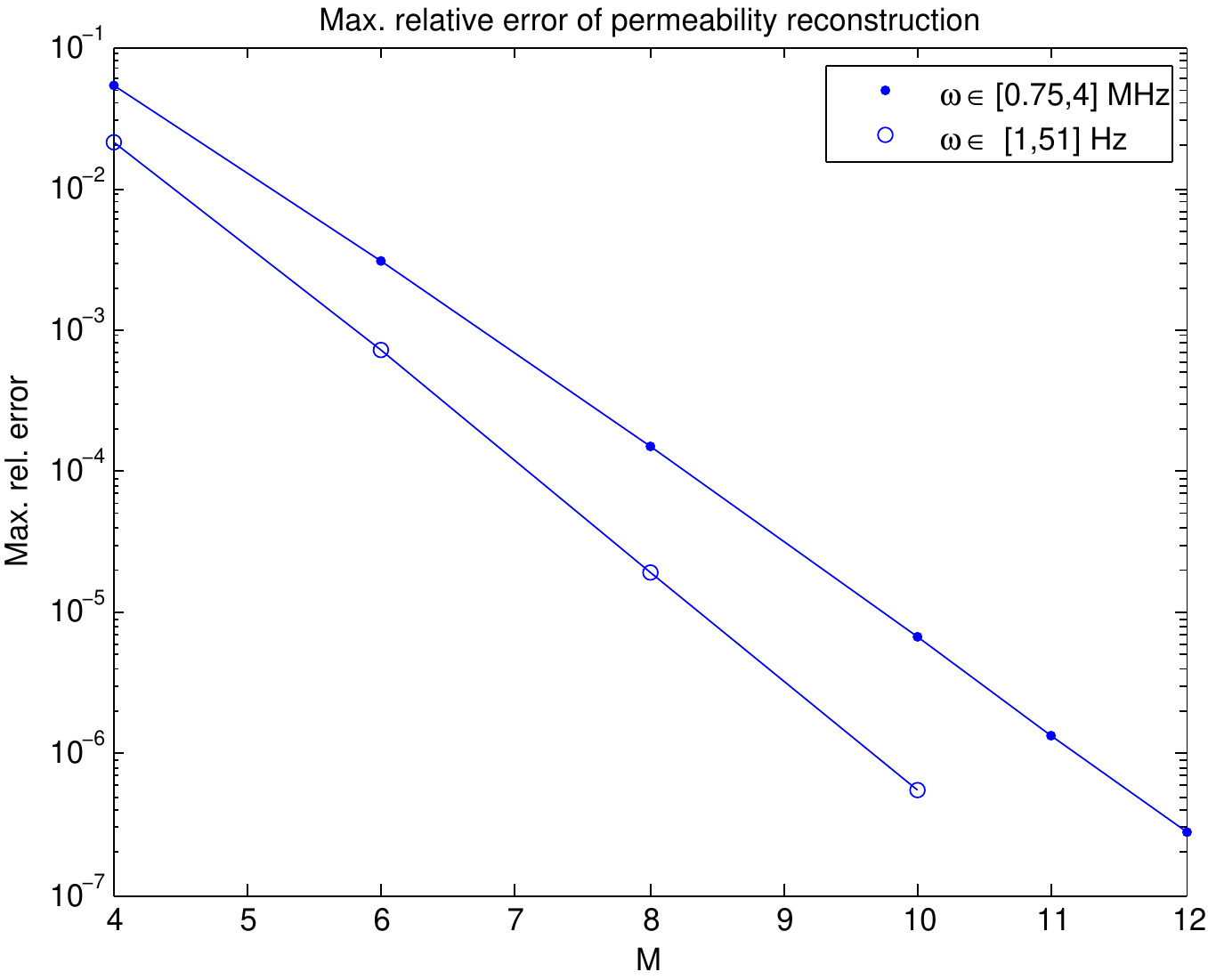}
\caption{\label{rel-err-1-51} $E_\infty$($P^M_{est}$) for $\omega\in [1, 51]$ Hz and $\omega\in [0.75, 4]$ MHz.}
\end{figure}
The second example is in the range from $0.75$ MHz to $4$ MHz, which is the spectrum range of the incident ultrasound wave used in \cite{fellah2004ultrasonic-wave} for studying cancellous bones. $E_\infty(P^M_{est})$ with respect to the exact $P$, which is $P^D$ is shown in Figure \ref{rel-err-1-51}. From data in this high frequency range, the moments are not well-approximated. However, it is interesting to note that the sum of residuals $\sum_{j=1}^{M'} r_j$, which serves as approximation to $\int_0^{\xi_p} \frac{d\lambda(u)}{u}$ are close to 1, as predicted by Theorem \ref{K-thm}. 

Examples 3 and 4 are taken from the spectral content of the incident waves used in \cite{lu2005wave-field-simu} for seismic wave modeling, where the memory term is handled by the shifted fractional derivative approach. In example 3, the frequency range is from 0 to 4 kHz while the range from 4 kHz to 180 kHz is considered in example 4. Due to the wide spreading of the frequency range, the reconstruction is not as efficient as the previous two cases. However, we note that a common feature in these two cases is that the increase in $E_\infty$ is due to the error in the frequency from $\omega_{min}$ to $\omega_{min}+\delta \omega$, i.e. between the first and the second sample points in the multipoint Pad\'{e} approximation. Improvement can be achieved by modifying the location of sample points. For example, in Figure \ref{0to4000-rel-err-M10-R+9p1}, the curve marked by red crosses is from equally spaced sample points for $M=10$, $\omega_{min}=40$ Hz, $\omega_{max}=4000$ Hz and $\delta \omega=440$ Hz, from which the value at $\omega=0$ can be approximated with relative error 1e-3; excluding the first interval, the maximum relative error drops from 6.25e-2 to 1.43e-4. 

\begin{figure}[h]
\centering
\includegraphics[width=0.95\textwidth]{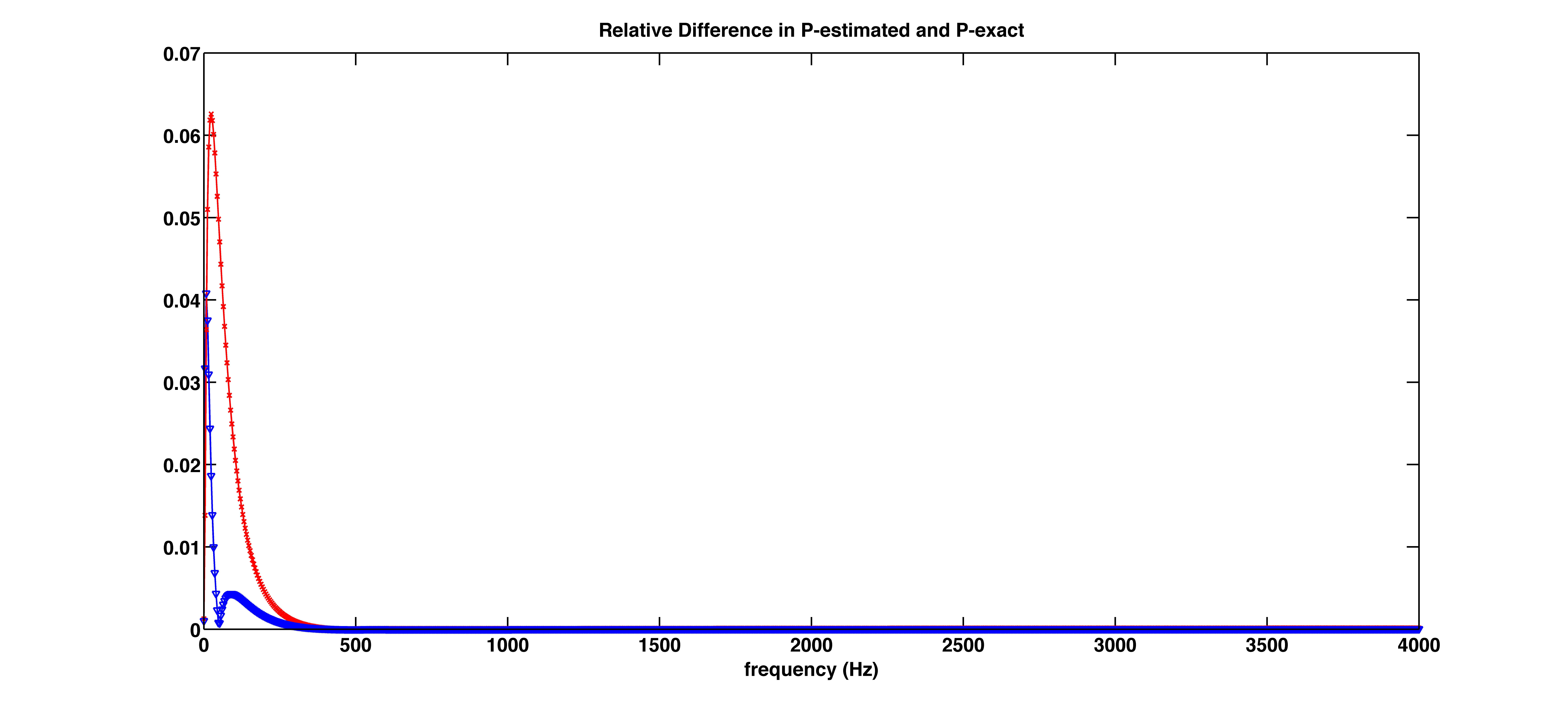}
\caption{\label{0to4000-rel-err-M10-R+9p1}Comparison of relative error of $P^{10}_{est}$ for different spacings of sample points. The red curve with crosses is from equally spaced sample points. The blue curve with triangles is from the modified approach.}
\end{figure}
%

Noting that the relative error peaks around $\delta \omega/10$ from the minimum frequency, we use 9 equally spaced sample points in the frequency range, which corresponds to $\delta \omega=500$, and an {\bf extra} sample point at $\omega_{min}+\delta \omega/10$.  The relative error from the modified approach is marked by blue triangles in Figure \ref{0to4000-rel-err-M10-R+9p1}. As indicated by (M) in Figure \ref{rel-err-0-4k} and Figure \ref{rel-err-4k-180k}, this modification brings down $E_\infty$($P^M_{est}$). The curve marked by (R) is from equally spaced sample points.
\begin{figure}[h]
\centering
\includegraphics[width=0.5\textwidth]{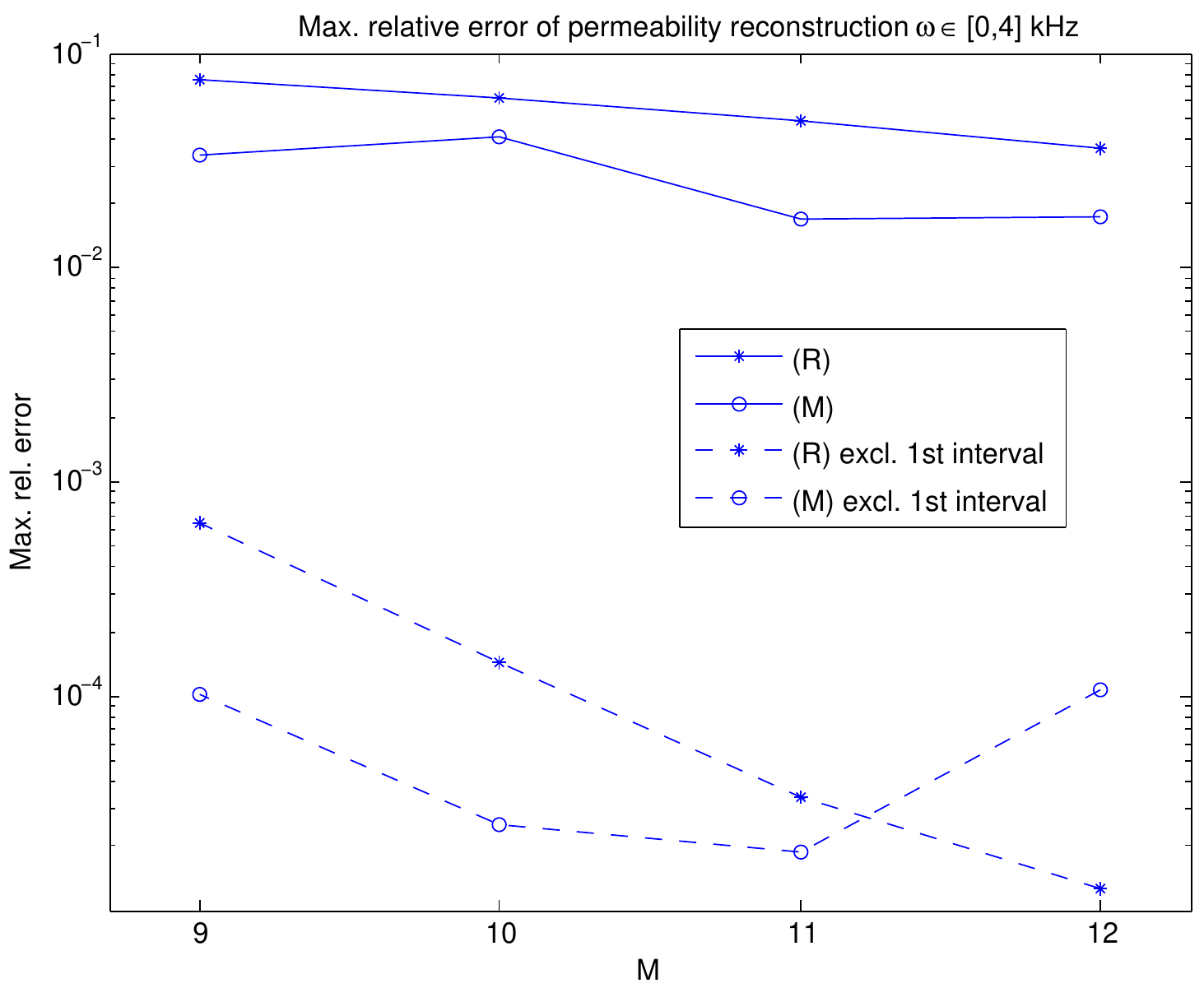}
\caption{\label{rel-err-0-4k} $E_\infty$($P^M_{est}$) for $\omega\in [0, 4]$ kHz. (R): equally spaced (M): Modified.}
\end{figure}
\begin{figure}[h]
\centering
\includegraphics[width=0.5\textwidth]{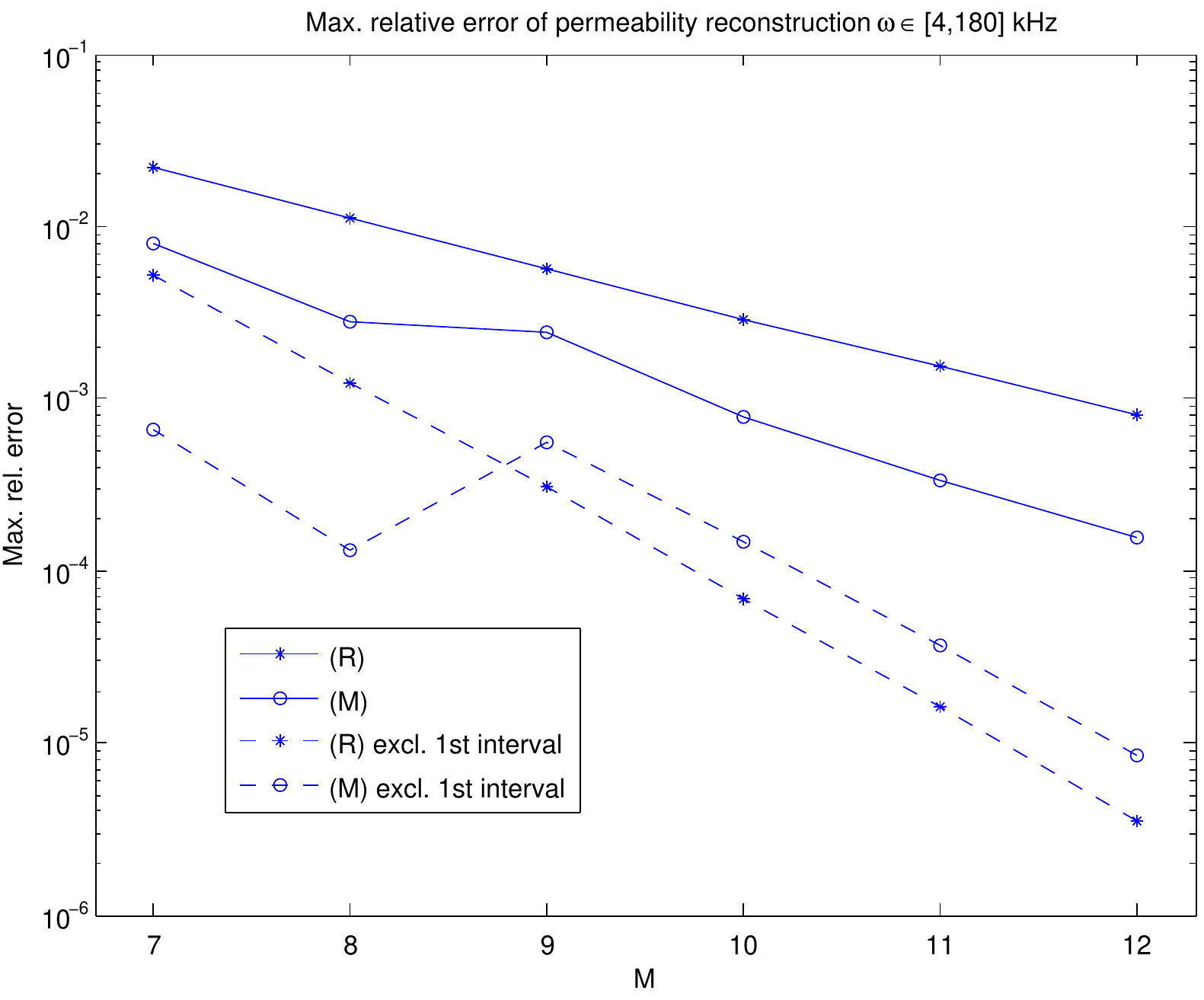}
\caption{\label{rel-err-4k-180k} $E_\infty$($P^M_{est}$) for $\omega\in [4, 180]$ kHz. (R): equally spaced (M): Modified.}
\end{figure}
\section{\label{recon-tortu}Reconstruction of dynamic tortuosity}
In the simulation of high-frequency wave propagation in poroelastic media, the time domain Darcy's laws that come from inverse Fourier transform of (\ref{Darcy-freq}) are part of the first-order formulation of balance law. To deal with the memory term, in \cite{carcione1996wave-propagatio} a phenomenological approach using generalized Zener kernels was proposed, yet not implemented, with relaxation times obtained by curve fitting. In this section, we show that the analytical structure of tortuosity $T$ in frequency domain can be utilized to calculate the parameters  needed in the dissipation kernels from a finite data set of $T(\omega_j)$, $j=1,\cdots, M$. 
\subsection{IRF for dynamic tortuosity $T(\omega)$}
We first note that (\ref{tortuosity-permeability}) implies
\beq
T(\omega)=-\left( \frac{\eta \phi F}{\rho_f \nu}\right)\left(\frac{1}{R(\xi)}\right)=-\frac{\alpha_\infty}{R(\xi)} 
\label{Tort-R}
\eeq
To derive the IRF, we embed $R(\xi)$ into a larger class $\mathcal{S}$.
\begin{definition} \cite{berg1980quelques-remarq}
A function $g: (0,\infty) \goto [0,\infty)$ is in $\mathcal{S}$ if it can be represented as follows
\[
g(\zeta)=a+\int_0^\infty \frac{d\sigma(t)}{\zeta+t},\, \zeta>0
\]
with constant $a\ge 0$ and a positive measure $d\sigma$ on $[0,\infty)$.
\end{definition}
Let $\zeta:=-\xi$ and $g(\zeta):=-R(\xi)$, which is well-defined for all $\zeta>0$ because all the singularities of $R(\xi)$ are confined in $\xi\in[0,\Theta_1]$. Hence $g(\zeta)\in \mathcal{S}$ with $a=0$ and $d\sigma=d\lambda$. It is known \cite{berg1980quelques-remarq} that $\frac{1}{\zeta g(\zeta) }\in \mathcal{S}$ if $g(\zeta)\in \mathcal{S}$, therefore
\[
\frac{\alpha_\infty}{\zeta g(\zeta)}=\frac{\alpha_\infty}{\xi R(\xi)}=a+\int_0^\infty \frac{d \sigma(t)}{-\xi+t} \mbox{ for some $a\ge 0$ and positive measure $d\sigma$}, \zeta>0
\]
Since every function in $\mathcal{S}$ can be analytically extended to the cut complex plane and $R(\xi)$ is analytic in $\field{C}\setminus [0,\Theta_1]$, the expression above is valid in $\field{C}\setminus [0,\Theta_1]$ and we conclude from (\ref{Tort-R}) that the tortuosity function $T(\omega)$ has the following IRF:
\begin{theorem}
The dynamic tortuosity $T(\omega)$ has the following integral representation formula for $\omega$ such that $-\frac{i}{\omega}\in \field{C}\setminus[0,\Theta_1]$
\beq 
T(\omega)=a\left(\frac{i}{\omega}\right) +\int_0^{\Theta_1}\frac{d\sigma(t)}{1-i\omega t}
\label{T-IRF}
\eeq
for some constant $a\ge 0$ and positive measure $d\sigma$.
\end{theorem} 
According to (\ref{tortuosity-permeability}), $T(\omega)$ has a pole at $\omega=0$ with strength $\frac{i\eta\phi}{\rho_f K_0}$, hence the $a$ in (\ref{T-IRF}) is
\[
a=\frac{\eta\phi}{\rho_f K_0} =\frac{\alpha_\infty}{C_2}
\]
Furthermore, $T(\omega)\goto \alpha_\infty$ as $\omega\goto \infty$, so $d\sigma$ has a Dirac mass at $t=0$ with strength $\alpha_\infty$. It is also interesting to note that (\ref{T-IRF}) implies
\beq
\mu_0(d\sigma)=\frac{\nu^2 \phi}{K_0^2 F} \mu_1(d\lambda)=\frac{\phi F}{\mu_0^2(d\lambda)} \mu_1(d\lambda)=\frac{\alpha_\infty}{\mu_0^2(d\lambda)} \mu_1(d\lambda)
\label{mu-relation}
\eeq
which can be easily seen by taking the limit $\lim_{\omega\goto 0} T(\omega)-\frac{ia}{\omega}$.

Suppose we have the data of $R(\xi_j)$, $\xi_j:=-\frac{i}{\omega_j}$, $j=1,\cdots,M$ at different nonzero frequencies $\omega_j\in\field{R}$, which can come from measurements of $K(\omega_j)$ or $T(\omega_j)$.  To reconstruct $T(\omega)$, we first recognize that  
\[
h(\xi):=a-\frac{\alpha_\infty}{\xi R(\xi)}=a+i\omega T(\omega)\]
 is a Stieltjes function with the symmetry $h(-\xi)=\overline{h(\xi)}$ for $\xi=-\frac{i}{\omega}$ and $\omega\in \field{R}$, i.e. measurement at $M$ different non-zero frequencies provide $2M$ data. Similar with Section \ref{recon-perm}, we can use these $2M$ data points to reconstruct $h(\xi)$ by using multipoint Pad\'{e} approximates
\beq
h(\xi_j) = \int_0^{\Theta_1} \frac{d\sigma(\Theta)}{\xi_j -\Theta} \approx [M-1/M]_h(\xi_j):=\frac{a_0+a_1 \xi_j+\cdots+a_{M-1} \xi_j^{M-1}}{1+b_1 \xi_j+\cdots+b_{M} \xi_j^{M}}, \,\, j=1,2,\cdots, 2M
\label{LM-Pade-h}
\eeq
Once the $[M-1/M]_h(\xi)$ is known, its partial fraction decomposition can be numerically obtained
\beq
[M-1/M]_h(\xi) = \sum_{j=1}^M \frac{r_j}{\xi-p_j},\, r_j>0, 0<p_j <\Theta_1
 \label{pole-res-def}
 \eeq  
and the dynamic tortuosity $T$ can be approximated in terms of residues $r_j$ and poles $p_j$, $j=1,\cdots,M$
 \beq
T(\omega) \approx  \frac{a}{-i\omega}   +  \sum_{j=1}^M \frac{r_j/p_j}{-i\omega +1/p_j}
\label{pol-res-T}
\eeq
Therefore, the tortuosity in time domain can be approximated as 
\beq
\mathcal{F}^{-1}[T](t)= a \delta(t)+\sum_{j=1}^M \left(\frac{r_j}{p_j}\right) e^{-\frac{t}{p_j}}, \, t\ge 0.
\label{T-quad-xi}
\eeq
with $p_j$ and $r_j$ defined in (\ref{pole-res-def}).
\subsection{Formulation and Algorithm for $T(\omega)$}
Given the permeability data $P(s_j)$, $s_j:=-i\omega_j$, $j=1,\cdots,M$ at different nonzero frequencies $\omega_j\in\field{R}$, we compute the data points for $D(s)$ defined as  
\beq
D(s):=T(\omega)-\frac{i a}{\omega}=\frac{\alpha_\infty}{sP(s)}-\frac{a}{s}=\int_0^{\Theta_1} \frac{d\sigma(\Theta)}{1+s\Theta}
\label{D-def}
\eeq
Note that 
\[
\lim_{s\goto 0} D(s)=\int_0^{\Theta_1} d\sigma(\Theta)<\infty.
\]
Using symmetry, there are $2M$ data points for reconstructing $D(s)$ by multipoint Pad\'{e} approximates
\beq
D(s_j) = \int_0^{\Theta_1} \frac{d\sigma(\Theta)}{1+s_j \Theta} \approx [M-1/M]_D(s):=\frac{a_0+a_1 s_j+\cdots+a_{M-1} s_j^{M-1}}{1+b_1 s_j+\cdots+b_{M} s_j^{M}}, \,\, j=1,2,\cdots, 2M
\label{LM-Pade-D}
\eeq
The linear system of $a_0,\cdots, a_M$ and $b_1,\cdots,b_M$ to be solved has the same structure as that in Section \ref{R-reconstruct} except $P(s_j)$ and $P_j$ in (\ref{linear-sys}) should be replaced with $D(s_j)$ and $D_j$, respectively. The numerical scheme is identical to the 4-step process described in Section \ref{R-reconstruct}.    

Once the $[M-1/M]_D(s)$ is known, its partial fraction decomposition can be numerically obtained
\beq
[M-1/M]_D(s) = \sum_{j=1}^M \frac{r_j}{s-p_j},\, r_j>0, p_j \in (-\infty, -\frac{1}{\Theta_1})
\label{res-pol-D}
\eeq 
and the dynamic tortuosity $T$ can be approximated in terms of the residues and poles of $[M-1/M]_D(s)$
\[
T(\omega) \approx  \frac{a}{-i\omega}   +  \sum_{j=1}^M \frac{r_j}{-i\omega -p_j}
\] 
Therefore, the tortuosity in time domain can be expressed as 
\beq
\mathcal{F}^{-1}[T](t)= a \delta(t)+\sum_{j=1}^M r_j  e^{p_j t}
\label{T-quad}
\eeq
where $r_j>0$ and $p_j<0$ are defined in (\ref{res-pol-D}) and $\delta(t)$ is the Dirac function.
\vspace{0.2in}

We use the JKD tortuosity function to demonstrate the idea. 
\subsection{Numerical Results for JKD tortuosity $T^D(\omega)$}
The function corresponding to $T^D$ via (\ref{D-def}) is 
\beq
D^D(s):= \frac{\alpha_\infty}{sP^D(s)}-\frac{a}{s}=\alpha_\infty+\frac{\alpha_\infty(\sqrt{1+C_1s}-1)}{C_2 s}=\int_0^{\Theta_1} \frac{d\sigma^D(\Theta)}{1+s\Theta}
\label{DD-def}
\eeq
We use the values of parameters in Table \ref{table-para} for the simulations. For these parameters, $a=\frac{\eta \phi}{\rho_f K_0}$=2.708895e03. 

Suppose $M'$ poles are retained after the algorithm, $M' \le M $, we reindex them and the corresponding residues to $\{(p_j,r_j)\}|_{j=1}^{M'}$.  The function $D(s)$ is then approximated by
\beq
D(s)\approx D^M_{est}(s):=\sum_{j=1}^{M'} \frac{r_j}{s-p_j}=\sum_{j=1}^{M'} \frac{r_j}{s-p_j}=\sum_{j=1}^{M'} \frac{-r_j/p_j}{1+s(-1/p_j)}
\label{D-est}
\eeq
and the moments $\mu_k(d\sigma^D)$, $k=0,1,2,\cdots$ by
\beq
\mu_k(d\sigma^D) \approx  (-1)^{k+1} \sum_{j=1}^{M'}  \frac{r_j}{(p_j)^{k+1}}.
\eeq
In terms of the poles and residues in (\ref{D-est}), the time domain tortuosity $T$ can be approximated as
\[
\mathcal{F}^{-1}\{T\}(t)\approx \frac{\eta \phi}{\rho_f K_0}\delta(t)+ \sum_{j=1}^{M'} r_j e^{p_j t}, \,\,r_j>0,\, p_j<0.
\]
We consider the same frequency ranges as in Section \ref{NuRe-K}. Table \ref{DT-moments} shows the reconstructed moments of $d\sigma^D$ from data in the frequency range from 1 to 51 Hz with multipoint Pad\'{e} approximants of various order $M$; the plot is demonstrated in Figure \ref{fig-table8}.  The exact values of moments are computed by first observing that  $D^D$ has a removable singularity at $s=0$ and its Taylor expansion near $s=0$ can be explicitly expressed as
\[
D^D(s)= \sum_{k=0}^\infty c_k s^{k}\]
with
\[
c_0= \alpha_\infty+ \frac{2K_0 \alpha_\infty^2}{\phi \Lambda^2}  \mbox{   and   }  c_k=\frac{(-1)^{k} C_1^{k+1} \alpha_\infty}{(k+1)! 2^{k+1} C_2} \prod_{j=1}^{k}(2j-1), \,\,k\ge 1
\]
Differentiating the IRF in (\ref{DD-def}) with respect to $s$ and compare both sides, we can express the moments of $d\sigma^D$ in terms of the Taylor coefficients of $D^D$ near $s=0$
\beq
\mu_k(d\sigma^D)=(-1)^k c_{k},\,\, k=0,1,2,\cdots
\label{sig-moment}
\eeq
%
\begin{table}[htb]
\begin{center}
\caption{\label{DT-moments} Moments of $d\sigma^D$ constructed from data in $[1, 51]$Hz}
{\small
\begin{tabular}{|c|c|c|c|c|c|}
\hline
& Exact & $M$=9 &  $M$=7 & $M$=5 &  $M$=3  \\
\hline
$M'$& & 8&6&4&2\\
\hline
$\mu_0$ & 0.2448054 e3 & 0.2448050 e3 & 0.2447906 e3 & 0.2442974 e3  & 0.2313087 e3\\
\hline
$\mu_1$ & 0.1096426 e2 & 0.1096420 e2 & 0.1096153 e0 & 0.1088098 e2 & 0.9286772 e1\\
\hline
$\mu_2$ & 0.9864792 e0 & 0.9864178 e0 & 0.9844351 e0 & 0.9460585 e0 & 0.5683834 e0\\
\hline
$\mu_3$ & 0.1109447 e0 & 0.1109081 e0 & 0.1100577 e0 & 0.9890620 e-1 & 0.3758564 e-1\\
\hline
$\mu_4$ & 0.1397472 e-1 & 0.1395836 e-1 & 0.1367750 e-1 & 0.1108704 e-1 & 0.2511712 e-2\\
\hline
$\mu_5$ & 0.1886006 e-2 & 0.1879998 e-2 & 0.1801971 e-2 & 0.1278614 e-2 & 0.1680769 e-3\\
\hline
$\mu_6$ & 0.2666530 e-3 & 0.2647426 e-3 & 0.2455669 e-3 & 0.1491975 e-3 & 0.1124921 e-4\\
\hline
$\mu_7$ & 0.3898596 e-4 & 0.3844154 e-4 & 0.3413699 e-4 & 0.1749433 e-4 & 0.7529144 e-6\\
\hline
$\mu_8$ & 0.5846091 e-5 & 0.5703658 e-5 & 0.4801441 e-5 & 0.2055462 e-5 & 0.5039304 e-7\\
\hline
$\mu_9$ & 0.8941761 e-6 & 0.8593632 e-6 & 0.6800106 e-6  & 0.2417042 e-6 & 0.3372839 e-8\\
\hline
\end{tabular}
}
\end{center}
\end{table}
\begin{figure}[h]
\centering
\includegraphics[width=0.5\textwidth]{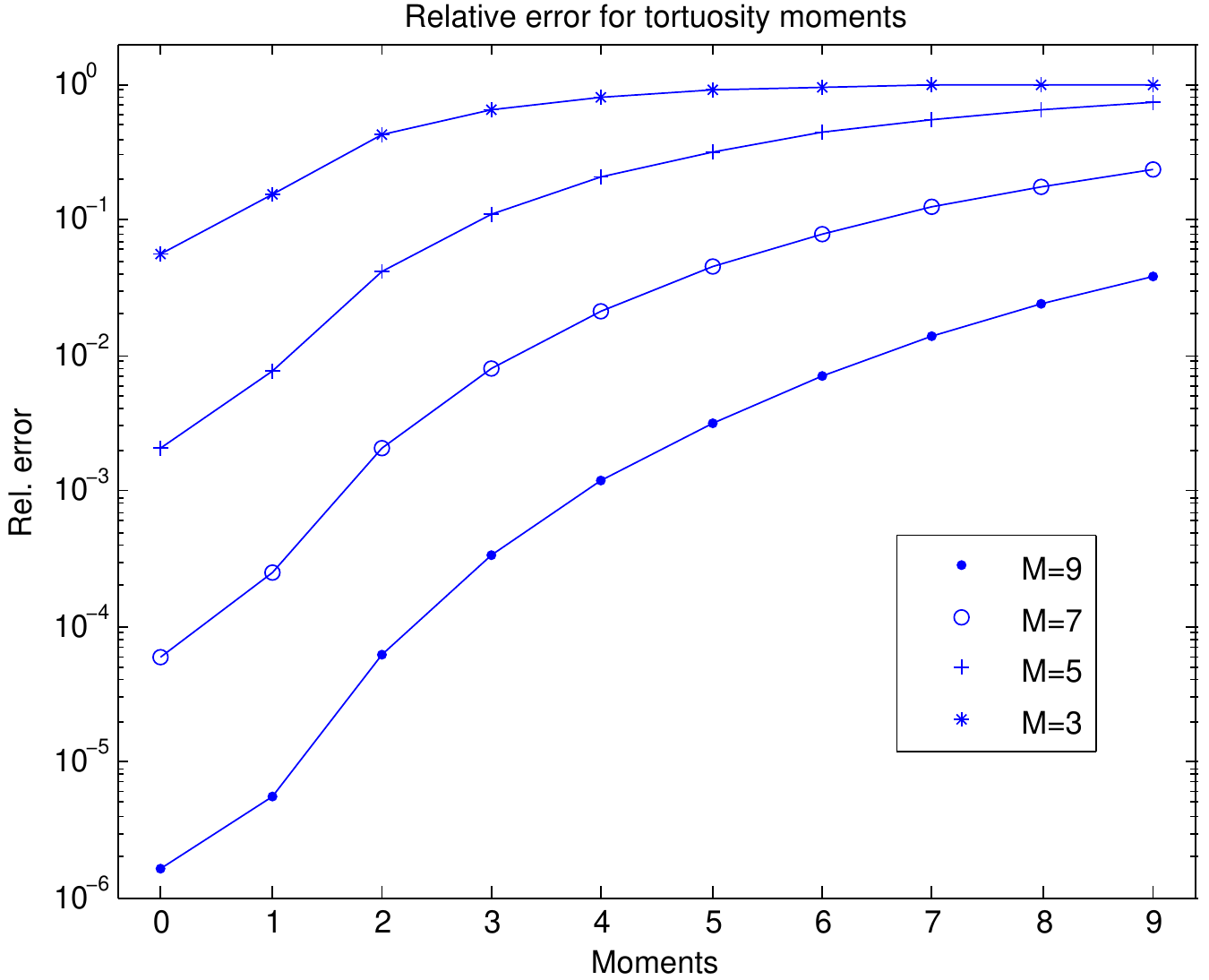}
\caption{\label{fig-table8}. Relative error of $\mu_1(d\sigma^D), \cdots, \mu_9(d\sigma^D)$ approximated from $M$ data points }
\end{figure}

Due to the important roles played by the poles and residues in handling the memory terms in Biot-JKD equations, we list $\{p_j,r_j\}_{j=1}^{M'}$ for $M=9$ and $M=7$ in Table \ref{1-51-pole-res} with $j=1,\cdots,M'$.
%
\begin{table}[htb]
\begin{center}
\caption{\label{1-51-pole-res} Poles and residues of $[M-1/M]_{D^D}(s)$ constructed from data in $[1, 51]$Hz, }
{\small
\begin{tabular}{|c|c|c|c|c|}
\hline
& \multicolumn{2}{|c|}{$M=9$} &   \multicolumn{2}{|c|}{$M=7$} \\
\hline
&$p_j$ &$r_j$ & $p_j$& $r_j$\\
\hline
$j=1$ & -1.706303 e3   &    6.240356 e4 & -9.475742 e2 & 4.524730 e4\\
\hline
$j=2$ & -1.988249 e2  &  6.285370 e3 & -1.183222 e2 & 4.649229 e3\\
\hline
$j=3$ &   -7.810948 e1  &     2.285702 e3  &  -4.647002 e1 & 1.807574 e3\\
\hline
$j=4$ & -4.120607 e1  &   1.252540 e3 & -2.284838 e1 & 1.018569 e3\\
\hline
$j=5$ & -2.398527 e1  &    8.133531 e1  &   -1.199244 e1 & 5.693354 e2\\
\hline
$j=6$ &  -1.450089 e1  &   5.338794 e1  &  -6.992178 e0 & 1.873584 e2 \\
\hline
$j=7$ & -9.133035 e0  &    2.945426 e1 & &\\
\hline
$j=8$ & -6.391875 e0  &    8.937405 e0 && \\
\hline
\end{tabular}
}
\end{center}
\end{table}
The maximum relative error of $D^M_{est}$ is defined as 
\[
E_\infty(D^M_{est}):=\max_{s \in [-i\omega_{min},-i\omega_{max}]} \frac{| D(s)-D^M_{est}| }{|D(s)|}
\]
and evaluated the same way as that of $P^M_{set}$ in Section \ref{NuRe-K}. The maximum relative error for frequency range $[1, 51]$Hz and $[0.75, 4]$ MHz, $[0, 4000]$Hz and $[4,180]$ kHz are listed in Figures \ref{rel-err-1-51-D}, \ref{rel-err-0-4k-D} and \ref{rel-err-4k-180k-D}, respectively. 
\begin{figure}[h]
\centering
\includegraphics[width=0.5\textwidth]{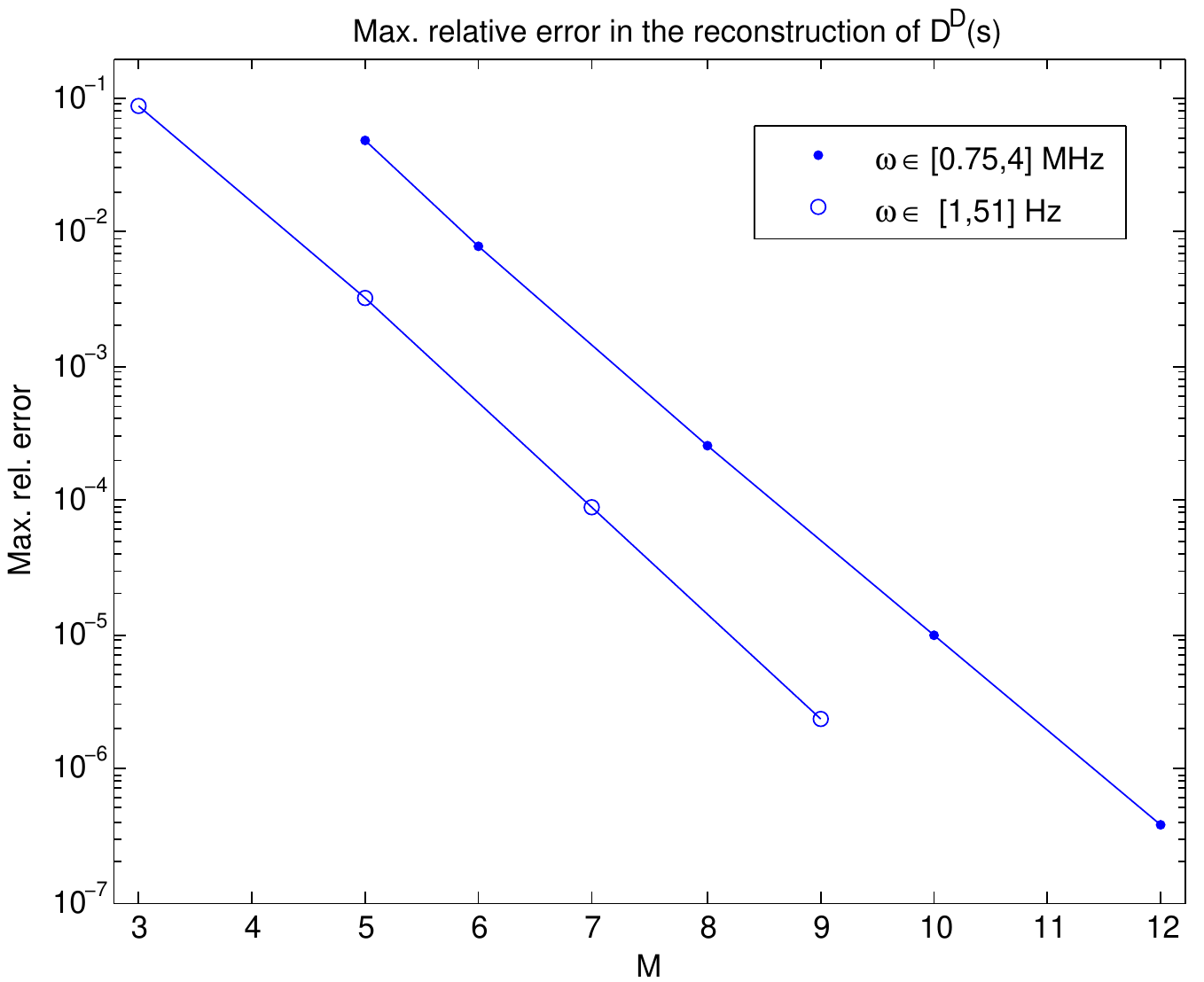}
\caption{\label{rel-err-1-51-D}. $E_\infty$($D^M_{est}$) for $\omega\in [1, 51]$ Hz and $\omega\in [0.75, 4]$ MHz }
\end{figure}
\begin{figure}[h]
\centering
\includegraphics[width=0.5\textwidth]{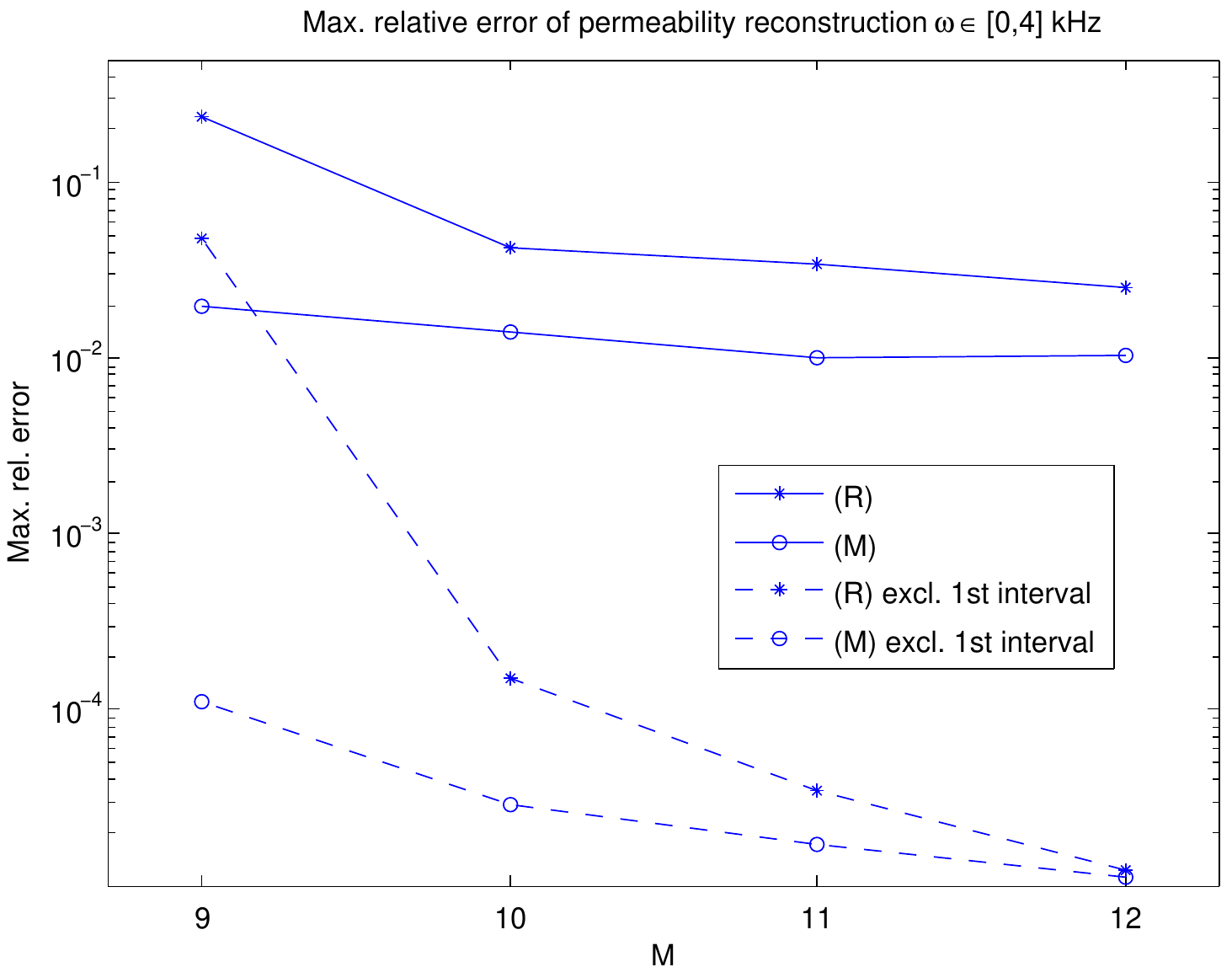}
\caption{\label{rel-err-0-4k-D}. $E_\infty$($D^M_{est}$) for $\omega\in [0, 4]$ kHz. (R): equally spaced (M): Modified.}
\end{figure}
\begin{figure}[h]
\centering
\includegraphics[width=0.5\textwidth]{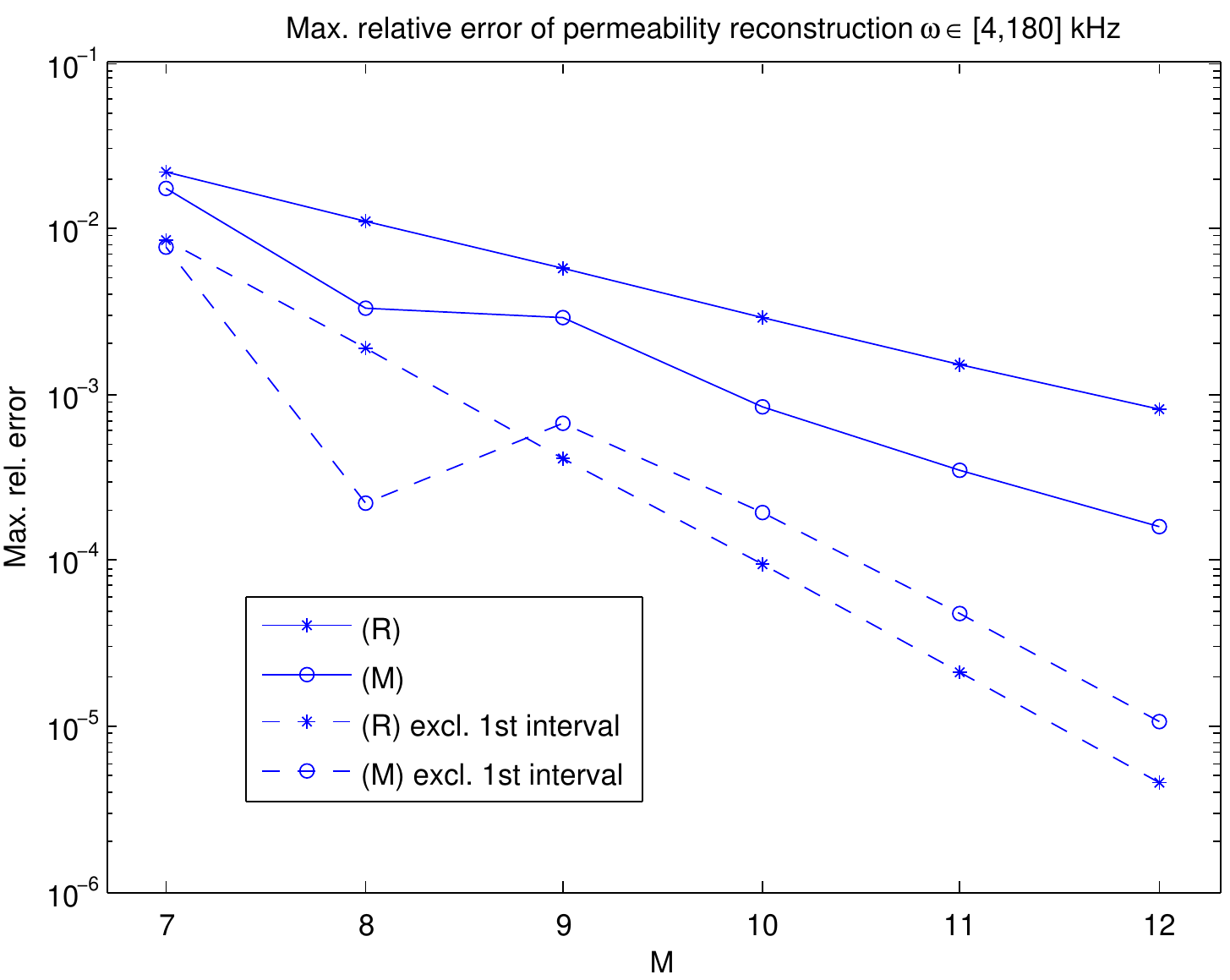}
\caption{\label{rel-err-4k-180k-D}. $E_\infty$($D^M_{est}$) for $\omega\in [4, 180]$ kHz. (R): equally spaced (M): Modified.}
\end{figure}
\section{\label{moment-remark}Remarks on the relation between moments and various effective parameters}
Using the IRFs, we can derive several relations between the moments of $d\lambda$, $d\sigma$ and various combinations of dynamic effective parameters. Recall that $d\lambda$ corresponds to the permeability function and $d\sigma$ to the tortuosity function. Since the dynamic permeability and dynamic tortuosity depend on both the frequency and the pore space geometry, the fact that the integrands in the IRFs are only functions of frequency implies that all the geometrical information must be encoded in the measures. The analyticity of both functions at $s=0$ enable the calculation of moments in terms of the coefficients of Taylor expansions there. For example, (\ref{P-function}) implies
\beq
\mu_0(d\lambda)=\frac{FK_0}{\nu}=\frac{\alpha_\infty K_0}{\phi \nu}
\label{zero-moment}
\eeq
The relation between the permeability and the tortuosity leads to the following relation between the moments of $d\lambda$ and $d\sigma$
\beq
\mu_p(d\lambda)=\frac{\mu_0(d\lambda)}{\alpha_\infty} \sum_{k+j=p-1} \mu_k(d\sigma)\,\mu_j(d\lambda),\, p\ge 1
\label{moment-relation}
\eeq
From (\ref{moment-relation}), the infinite-frequency $\alpha_\infty$ can be expressed in terms of moments
\beq
\alpha_\infty=\frac{\mu_0(d\sigma) \mu_0^2(d\lambda)}{\mu_1(d\lambda)}
\label{alpha-moment}
\eeq
For the JKd permeability and tortuosity, the measure $d\sigma^D$ satisfies 
\beq
\mu_0(d\sigma^D)=\int_0^{\Theta_1} d\sigma^D(\Theta) = \lim_{s\goto 0} D^D(s)=\alpha_\infty + \frac{2K_0 \alpha_\infty^2}{\phi \Lambda^2}
\label{mu0-sig}
\eeq
Suppose we can reconstruct $\mu_0(d\lambda^D)$ and $\mu_0(d\sigma^D)$ from low frequency data of permeability $K(\omega)$ and that we know the porosity $\phi$ and pore fluid kinetic viscosity $\nu$, then we can recover $\alpha_\infty$ from the fact that $\mu_0(d\lambda^D)=\frac{\alpha_\infty K_0}{\phi \nu}$ because $K_0$ can be obtained easily through extrapolation of low frequency data. Once $\alpha_\infty$ is known,  (\ref{mu0-sig}) can be used to recover $\Lambda$, which is a weighted pore volume-to-surface ratio that provides a measure of the dynamically connected part of the pore region \cite{avellaneda1991rigorous-link-b}. We note that (\ref{mu0-sig}) implies the following for the JKD model
\beq
\Lambda=\sqrt{\frac{2K_0 \alpha_\infty^2}{\phi[\mu_0({d\sigma^D})-\alpha_\infty]}} = \sqrt{\frac{2K_0 \alpha_\infty}{\phi[\frac{\mu_1(d\lambda^D)}{\mu_0^2(d\lambda^D)}-1]}} 
\label{JKD-Lambda}
\eeq

Formulas (\ref{zero-moment}), (\ref{alpha-moment}) and (\ref{JKD-Lambda}) show exactly how microstructural information affects the effective parameters $\alpha_\infty$ and $\Lambda$ through moments. 
\section{\label{conclusion}Conclusion}
In this paper, we derived the integral representation formula (IRF) for dynamic tortuosity $T(\omega)$ in general form; we show that $T(\omega)$ can be written as the sum of a function with a simple pole at 0 and a Stieltjes function. Utilizing the analytic structure of this IRF and the IRF of permeability $K(\omega)$ derived in \cite{avellaneda1991rigorous-link-b}, an algorithm based on multipoint Pad\'{e} approximation of Stieltjes functions is proposed for constructing  $K(\omega)$ and  $T(\omega)$ from the values of permeability at district frequencies. Taking into account the symmetry of Stieltjes functions, only $M$ different frequencies, instead of $2M$, are needed for constructing the $[M-1/M]$ approximant. It is demonstrated that the moments of both the measures in the IRFs of $K(\omega)$ and $T(\omega)$ can be estimated to high accuracy from low frequency data. The capability of this algorithm for recovering the moments can be utilized to compute the inf-tortuosity $\alpha_\infty$ through (\ref{zero-moment}) using the low-frequency permeability data, if the viscosity of pore fluid $\nu$ is known because $K_0$ can be approximated very well from $P^M_{est}(s)$. 
\vspace{0.1in}

Furthermore, if the JKD model is used, the microstructure-dependent parameter $\Lambda$ can be recovered by the formula in (\ref{JKD-Lambda}). It is interesting to note that the empirical formula for $\Lambda$ suggested by JKD \cite{Pride1992Modeling-The-Dr} is  $\Lambda\approx\sqrt{\frac{2\alpha_\infty K_0}{\phi/4}}$. Comparing this with (\ref{JKD-Lambda}), which is exact, this empirical formula corresponds to the assumption that $\frac{\mu_1(d\lambda^D)}{\mu_0^2(d\lambda^D)}=\frac{5}{4}$, which is not always true and obviously not satisfied by the moments calculated in this paper.  
\vspace{0.1in}

We have also shown that the JKD permeability $K^D(\omega)$ can indeed be represented by a probability measure in its IRF, as is predicted by the general result in \cite{avellaneda1991rigorous-link-b}.
\vspace{0.1in}

The results of numerical experiments conducted on the frequency ranges taken from  the literature in biomechanics for bone \cite{sebaa2006ultrasonic-char} and seismology \cite{lu2005wave-field-simu} are presented. For the bandwidth spreading less than two others of magnitude, the proposed algorithm with equally spacing interpolating points achieve approximates with high accuracy. From the last two numerical examples, we see that the approximation is of good accuracy away from the first interval; the max. relative error can be greatly improved by adding one sample point close to the lowest frequency to the equally spaced points. This implies that the location of sample points play an important role in the approximation and will be the topic of future investigation. Another way to handle wide frequency range can be to divide it into shorter intervals and do local approximation. The advantage of our reconstruction scheme is two-fold. First of all, it provides high accuracy interpolation of the permeability/tortuosity data without assuming anything beyond the fact that they are related to Stieltjes functions and hence is more general than the JKD model, which assumes specific forms of the dynamic tortuosity functions. Secondly, the time domain representation such as  (\ref{T-quad}) provides an efficient way for numerically handling the memory terms that appears in the time domain numerical simulation for wave propagation in poroelastic materials.
\vspace{0.1in}  
  
\noindent
{\bf Acknowledgement:} This research is partially sponsored by ARRA-NSF-DMS Math. Biology Grant 0920852.

%
\end{document}